\documentclass[11pt,letterpaper,reqno]{amsart}
\usepackage{amssymb}
\usepackage{esint}
\usepackage{mathtools}
\usepackage{fullpage}
\usepackage{hyperref}
\usepackage{enumitem}
\usepackage{setspace}
\usepackage{caption,subcaption,overpic}
\usepackage[normalem]{ulem}

\usepackage[english]{babel}
\usepackage[utf8]{inputenc}

\usepackage[usenames,dvipsnames]{color}
\usepackage[all]{xy} % for complicated commutative diagrams
\usepackage[alphabetic,backrefs,lite]{amsrefs}

\newtheorem{mainthm}{Theorem}
\newtheorem{mainconj}[mainthm]{Conjecture}

\theoremstyle{plain}
\newtheorem{theorem}{Theorem}%[section]
\newtheorem{lemma}[theorem]{Lemma}

\newtheorem{proposition}[theorem]{Proposition}

\theoremstyle{definition}
\newtheorem{definition}[theorem]{Definition}
\newtheorem{notation}[theorem]{Notation}
\newtheorem{example}[theorem]{Example}
\theoremstyle{remark}
\newtheorem{remark}[theorem]{Remark}

\newtheorem{remarks}[theorem]{Remarks}

\newcommand{\R}{\mathbb{R}}
\newcommand{\Z}{\mathbb{Z}}
\newcommand{\C}{\mathbb{C}}

\newcommand{\Q}{\mathbb{Q}}
\newcommand{\F}{\mathbb{F}}
\newcommand{\Aff}{{\mathbb A}}
\newcommand{\PP}{{\mathbb P}}
\newcommand{\calO}{{\mathcal O}}
\newcommand{\Sphere}{\mathbb{S}}

\newcommand{\bS}{\mathbf{S}}

\newcommand{\nakeddata}{\mathcal{T}}
\newcommand{\decorateddata}{\mathcal{D}}

\DeclareMathOperator{\Spec}{Spec}
\DeclareMathOperator{\Proj}{Proj}
\newcommand{\into}{\hookrightarrow}
\newcommand{\defi}[1]{\textsf{#1}} % for defined terms
\newcommand{\isom}{\simeq}
\newcommand{\frakp}{{\mathfrak p}}

\newcommand{\vect}[1]{\mathbf{#1}}

\newcommand{\ip}[2]{\left\langle#1,#2\right\rangle}

\newcommand{\eps}{\varepsilon}
\renewcommand{\phi}{\varphi}
\newcommand{\abs}[1]{\left| #1 \right|}

\DeclareMathOperator{\Irr}{Irr}
\DeclareMathOperator{\Gal}{Gal}
\DeclareMathOperator{\Bl}{Bl}
\DeclareMathOperator{\Sym}{Sym}

\newcommand{\pairs}{\psi}
\newcommand{\stiffnesstensor}[1]{E(#1)}
\newcommand{\posstiffnesstensor}[1]{E^+(#1)}
\newcommand{\orthorhombic}{E_{\mathrm{orth}}(3)}
\newcommand{\posorthorhombic}{E_{\mathrm{orth}}^+(3)}
\newcommand{\monoclinic}{E_{\mathrm{mono}}(3)}
\newcommand{\posmonoclinic}{E_{\mathrm{mono}}^+(3)}
\newcommand{\spoly}[1]{P_{#1}}
\newcommand{\spolymap}{F}
\newcommand{\dummy}{\,{\cdot}\,}

\title{Reconstruction of anisotropic stiffness tensors from partial data around one polarization}
\author{Maarten V. de Hoop, Joonas Ilmavirta, \\ Matti Lassas, and Anthony V\'arilly-Alvarado}
\address{Simons Chair in Computational and Applied Mathematics and Earth Science, Rice University, 6100 S.\ Main St., Houston, TX 77005, USA}
\email{mvd2@rice.edu}
\address{Department of Mathematics and Statistics, University of Jyv\"askyl\"a, P.O.Box 35 (MaD) FI-40014 University of Jyv\"askyl\"a, Finland}
\email{joonas.ilmavirta@jyu.fi}
\address{Department of Mathematics and Statistics, University of Helsinki, P.O. Box 68 (Gustaf H\"allstr\"omin katu 2b) FI-00014, University of Helsinki, Finland}
\email{Matti.Lassas@helsinki.fi}
\address{Department of Mathematics MS 136, Rice University, 6100 S.\ Main St., Houston, TX 77005, USA}
\email{av15@rice.edu}

\makeatletter
\@namedef{subjclassname@2020}{\textup{2020} Mathematics Subject Classification}
\makeatother

\date{\today}

\begin{document}

\begin{abstract}
We study inverse problems in anisotropic elasticity using tools from algebraic geometry. The singularities of solutions to the elastic wave equation in dimension $n$ with an anisotropic stiffness tensor have propagation kinematics captured by so-called slowness surfaces, which are hypersurfaces in the cotangent bundle of $\R^n$ that turn out to be algebraic varieties. Leveraging the algebraic geometry of families of slowness surfaces we show that, for tensors in a dense open subset in a space of anisotropic two-dimensional stiffness tensors, a small amount of data around one polarization in an individual slowness surface uniquely determines the entire slowness surface and its stiffness tensor. In three dimensions, for generic orthorhombic and monoclinic stiffness tensors, a small number of anomalous companions give rise to the same slowness surface; nevertheless, we conjecture that in the most anisotropic setting (triclinic) the tensor is unique, as in two dimensions. The partial data needed to determine a tensor arises naturally from seismological measurements or geometrized versions of seismic inverse problems. Additionally, we explain how the reconstruction of the stiffness tensor can be carried out effectively, using Gr\"obner bases.  Our uniqueness or finiteness results fail for symmetric materials (e.g., fully isotropic), evidencing the counterintuitive claim that inverse problems in elasticity can become more tractable with increasing asymmetry. 
\end{abstract}

\maketitle

%%%%%%%%%%%%%%%%%%%%%%%%%%%%%%%%
\setcounter{section}{-1}
\section{Prelude}
In an elastic medium, characterized by a stiffness tensor field, an excitation will produce waves of different polarizations propagating at different speeds. We set out to reconstruct the stiffness tensor field in a body from arrival time boundary observations of the fastest waves passing through it. This is an intricate problem. In dimension three it requires unraveling a large number of tensor components from scalar measurements. It is interesting even when the stiffness tensor field is independent of position (homogeneous): In this setting, it is impossible to solve when the medium is rotation invariant (isotropic).

Surprisingly, anisotropy of a medium makes reconstruction possible: We show that the same travel time measurements \emph{uniquely} determine a generic fully anisotropic homogeneous stiffness tensor in two dimensions, and for a generic orthorhombic or monoclinic stiffness tensor in three dimensions, only a small and understood number of tensors are compatible with given measurements. This uniqueness or finiteness is inaccessible with tools from analysis and differential geometry alone, and builds instead on the non-trivial algebraic geometry of the so-called slowness surfaces. Most natural media are anisotropic, and the often applied simplifying assumption of isotropy leads to unnecessary modeling restrictions and non-unique solutions to inverse problems.

To understand an inhomogeneous medium, which is outside the scope of this paper, one would have to understand the global structure of a bundle of slowness surfaces (the characteristic varieties of the elastic wave equation). Our results pertain to the analysis of the fibers of this bundle, a crucial first step.

%%%%%%%%%%%%%%%%%%%%%%%%%%%%%%%%
\section{Introduction}

Inverse problems in anisotropic elasticity are notoriously challenging: a lack of natural symmetry leaves one with few tools to approach them.  In this paper we embrace and harness asymmetry with the help of algebraic geometry, and develop a method to address inverse problems around the reconstruction of anisotropic stiffness tensors from a relatively small amount of empirical data. Along the way we prove surprisingly strong uniqueness and finiteness results in anisotropic elastic inverse problems, aided by the specific properties of olivine and muscovite, common materials in the Earth's mantle and crust, respectively. We view our results as the beginning of a fruitful interaction between the fields of inverse problems and modern algebraic geometry.

Microlocal analysis, describing the geometry of wave propagation, and algebraic geometry, describing the geometry of zero sets of polynomials, become linked through the \defi{slowness polynomial}, which is the determinant of the principal symbol of the elastic wave operator. The vanishing set of this polynomial is the \defi{slowness surface}, which describes the velocities of differently polarized waves in different directions. Notably, we show that for a \emph{generic} anisotropic material (in full generality in two dimensions and in specific symmetry types in three dimensions as described above)
\smallskip
\begin{enumerate}[leftmargin=*]
\item the polarizations of waves travelling through the material, corresponding to different sheets of the slowness surface, are coupled: a small Euclidean open subset of the slowness surface for a single polarization determines the whole slowness surface for all polarizations;
\smallskip
\item one can reconstruct the stiffness tensor field of the material from a slowness polynomial. 
\end{enumerate}

%%%%%%%%%%%
\subsection{The model}
\label{sec:stiffness-and-EWE}

We work in $\R^n$ for any~$n$; physical applications arise typically when $n=2$ or~$3$.

\subsubsection{Waves in anisotropic linear elasticity}

Linear elasticity posits that when a material is strained from a state of equilibrium, the force returning the system to equilibrium depends linearly on the displacement experienced. The displacement is described by the strain tensor (a symmetric $n\times n$ matrix) and the restoring force by a stress tensor (also a symmetric $n\times n$ matrix). Hooke's law, valid for small displacement to good accuracy, states that stress depends linearly on strain, and the coefficients of proportionality are gathered in a stiffness tensor. To be within the framework of linear elasticity, the displacement should be small, but one can also see the linear theory as a linearization of a more complicated underlying model.

Thus, the stiffness tensor of a material at a point $x \in \R^n$ is a linear map $\vect{c} = \vect{c}(x)\colon \R^{n\times n} \to \R^{n\times n}$ mapping \defi{strain} $\eps\in\R^{n\times n}$ (describing infinitesimal deformations) to \defi{stress} $\sigma\in\R^{n\times n}$ (describing the infinitesimal restoring force), given in components as
\begin{equation*}
\sigma_{ij} = \sum_{k,l=1}^n c_{ijkl} \eps_{kl}.
\end{equation*}
Since both $(\sigma_{ij})$ and $(\eps_{kl})$ are symmetric $n\times n$ matrices, we must have
\begin{equation}
\label{eq:minor-symmetry}
c_{ijkl} = c_{jikl} = c_{ijlk}.
\end{equation}
This is the so-called \defi{minor symmetry} of the stiffness tensor. In addition, the stiffness tensor is itself a symmetric linear map between symmetric matrices, which can be encoded in components as
\begin{equation}
\label{eq:major-symmetry}
c_{ijkl} = c_{klij}.
\end{equation}
This condition is known as the \defi{major symmetry} of the stiffness tensor.  Scaling by the density $\rho = \rho(x)$ of a material does not affect any of these properties, which leads to the \defi{reduced stiffness tensor} $\vect{a} = \vect{a}(x)$, whose components are $a_{ijkl} \coloneqq \rho^{-1}c_{ijkl}$.  Finally, a stiffness tensor is positive definite; combining the symmetries~\eqref{eq:minor-symmetry} and~\eqref{eq:major-symmetry} and positivity leads to the following definition formalizing the properties just described.

\begin{definition}
\label{def:stiffness-tensor}
We say that $\vect{a} = (a_{ijkl})\in\R^{n\times n\times n\times n}$ is a \defi{stiffness tensor} if 
\begin{equation}
\label{eq:symmetry}
a_{ijkl} = a_{jikl} = a_{klij}.
\end{equation}
If, in addition, for every non-zero symmetric matrix $A\in\R^{n\times n}$ we have
\begin{equation}
\label{eq:ellipticity}
\sum_{i,j,k,l=1}^n a_{ijkl}A_{ij}A_{kl} > 0,
\end{equation}
then we say that~$\vect{a}$ is \defi{positive}.
\end{definition}

\begin{notation}
\label{notation:stiffnesstensor}
The \defi{set of stiffness tensors} in $\R^{n\times n\times n\times n}$ is denoted $\stiffnesstensor{n}$, and the subset of positive ones is denoted $\posstiffnesstensor{n}$. Both sets carry a natural Euclidean topology.
\end{notation}

\begin{notation}
\label{nt:orthomono}
In dimension $n=3$ we name two special symmetry types of stiffness tensors important in materials science and mineralogy:
\smallskip
\begin{itemize}[leftmargin=*]
    \item A stiffness tensor $\vect{a}\in\stiffnesstensor{3}$ is called \defi{monoclinic} if
    \[a_{1123} = a_{1112} = a_{2223} = a_{2212} = a_{3323} = a_{3312} = a_{2313} = a_{1312} = 0.\]
    The space of monoclinic stiffness tensors is denoted by $\monoclinic\subset\stiffnesstensor{3}$. Examples of materials with monoclinic stiffness tensors include gypsum, borax, muscovite, and some synthetic ceramics and polymers.
    \smallskip
    \item A monoclinic stiffness tensor $\vect{a}$ is called \defi{orthorhombic} if additionally
    \[
    a_{1113} = a_{2213} = a_{3313} = a_{2312} = 0.
    \]
    The space of orthorhombic stiffness tensors is denoted by $\orthorhombic\subset\monoclinic$. Examples of materials with orthorhombic stiffness tensors include olivine, aragonite, and barite. Olivine is the primary component of the Earth's upper mantle.
\end{itemize}
\smallskip
The positive definite subsets of these vector spaces are respectively denoted $\posmonoclinic$ and $\posorthorhombic$.
\end{notation}

\begin{remark}
The symmetry classes of stiffness tensors can be understood in terms of the group of symmetries of the underlying Euclidean space (a subgroup of~$O(n)$) leaving the stiffness tensor invariant. The group always includes the two-element group~$G$ generated by the reflection $\vect{p}\mapsto-\vect{p}$. In the orthorhombic case the only non-trivial symmetry is $(x,y,z)\mapsto(x,-y,z)$, a flip across the $xz$-plane. In the monoclinic case flips across the other coordinate planes are also included. These symmetry classes therefore require the material's symmetry planes to align with the coordinate planes. This is not the case for triclinic stiffness tensors (minimal symmetry) and isotropic ones (maximal symmetry), corresponding to the two trivial subgroups of~$O(n)/G$.
\end{remark}

Denote the displacement from equilibrium of a material with stiffness tensor~$\vect{c}$ at point $x\in\R^n$ and time $t\in\R$ by $u(x,t)\in\R^n$.
The time evolution of $u(x,t)$ is governed by the \defi{elastic wave equation}
\begin{equation}
\label{eq:EWE}
\sum_{j,k,l=1}^n \frac{\partial}{\partial {x_j}} \bigg( c_{ijkl}(x) \frac{\partial}{\partial {x_k}} u_l(x,t) \bigg) - \rho(x)\frac{\partial^2}{\partial t^2}u_i(x,t) = 0,
\end{equation}
which can also be written as $\square u=0$, where $\square=\square_{c,\rho}$ is the matrix-valued elastic wave operator.

\subsubsection{Singularities and the principal symbol}

Let $\xi\in T^*_x\R^n$ be the momentum variable dual to~$x \in\R^n$ and let $\omega\in\R$ be the dual variable of time $t\in\R$. Following the terminology of microlocal analysis, a function $u(x,t)$ is said to be \defi{singular} at a point $(x_0,t_0)$ if $u(x,t)$ is not a $C^\infty$-smooth function in any neighborhood of the point $(x_0,t_0)$. A more precise description of singularities is given by the wave front set WF$(u)$ of the function $u(x,t)$, which consists of the points $(x,\xi,t,\omega)$ for which $u(x,t)$ is non-smooth at $(x,t)\in \R^n\times \R$ in the direction $(\xi,\tau)\in \R^n\times \R$. See \cite{Hormandervol1} for more details.

The singularities of $u(x,t)$ propagate by the null bicharacteristic flow of the matrix-valued principal symbol of~$\square$:
\begin{equation*}
\sigma(\square)_{il}(x,\xi,t,\omega) = -\sum_{j,k=1}^nc_{ijkl}(x)\xi_j\xi_k + \rho(x)\delta_{il}\omega^2.
\end{equation*}
A propagating singularity is annihilated by the principal symbol, so a point 
\[
(x,t,\xi,\omega)\in (\R^n\times \R)\times (\R^n\times \R)
\]
can be in the wave front set of a solution $u(x,t)$ of the elastic wave equation $\square u(x,t)=0$ only when
\begin{equation}
\label{eq:det(ps)=0}
\det( \sigma(\square)(x,\xi,t,\omega) ) = 0.
\end{equation}
Due to the homogeneity of the equation of motion, the frequency of oscillation has no effect on the propagation of singularities. It is therefore convenient to replace the momentum $\xi\in T^*_x\R^n$ with the \defi{slowness vector} $\vect{p}= (p_1,\dots,p_n) \coloneqq \omega^{-1}\xi\in T^*_x\R^n$.

To make~\eqref{eq:det(ps)=0} explicit, we recall that the \defi{Christoffel matrix} $\Gamma(x,\vect{p})$ is the $n\times n$ matrix whose $il$-th entry is
\begin{equation}
\label{eq:ChristoffelEntry}
\Gamma(x,\vect{p})_{il} = \sum_{j,k=1}^na_{ijkl}(x)p_jp_k.
\end{equation}
By~\eqref{eq:minor-symmetry} and~\eqref{eq:major-symmetry}, the Christoffel matrix is symmetric. If the stiffness tensor is positive, then the Christoffel matrix is positive definite. With this notation, the principal symbol becomes simply
\begin{equation*}
\sigma(\square) = - \omega^2\rho(x)(\Gamma(x,\vect{p})-I_n),
\end{equation*}
where~$I_n$ is the $n\times n$ identity matrix, and condition~\eqref{eq:det(ps)=0} can be rewritten as
\begin{equation}
\label{eq:det(G-I)=0}
\det( \Gamma(x,\vect{p})-I_n ) = 0.
\end{equation}
See Figure~\ref{fig:2d-slowness-example} for an example of the set of~$\vect{p}$'s that satisfy this condition.

\begin{remark}
Equation~\eqref{eq:det(G-I)=0} can also be argued physically by freezing the stiffness tensor~$c(x)$ and density~$\rho(x)$ to constant values~$c(x_0)$ and~$\rho(x_0)$ and writing a plane wave Ansatz for the displacement field. This is less rigorous but relies on the same underlying ideas and leads to the same condition. The principal symbol can also be understood as a description of how the operator acts on plane waves.
\end{remark}

\subsubsection{Polarization, slowness, and velocity}

\definecolor{color1}{rgb}{.6,.0,.6}
\definecolor{color2}{rgb}{.0,.8,.0}

\begin{figure}
\centering
\begin{subfigure}{.45\textwidth}
  \centering
  \includegraphics[width=.9\linewidth]{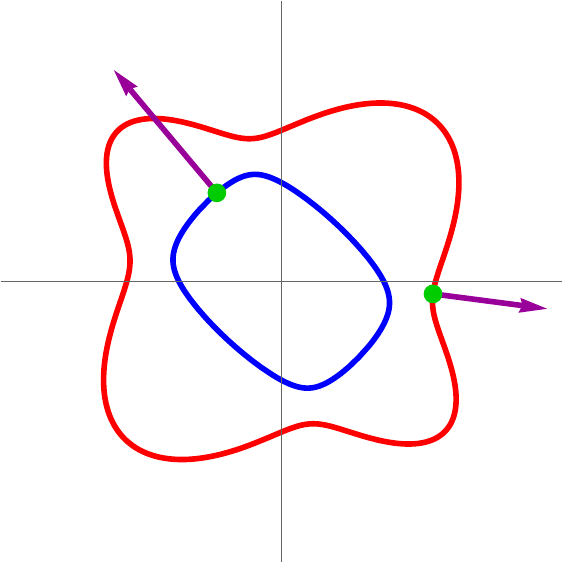}
  \caption{Slowness curve $S_x\subset T^*_x\R^2$, consisting of two colored branches.
  The Christoffel matrix $\Gamma(\vect{p})$ has two eigenvalues, $0<\lambda_1(\Gamma(\vect{p}))<\lambda_2(\Gamma(\vect{p}))$. The \textcolor{blue}{blue branch} corresponds to the values of $\vect{p}$ for which $\lambda_1(\Gamma(\vect{p}))=1$, and the \textcolor{red}{red branch} corresponds to $\lambda_2(\Gamma(\vect{p}))=1$.}
  \label{fig:sub1}
\end{subfigure}%
\hspace{1em}%
\begin{subfigure}{.45\textwidth}
  \centering
  \includegraphics[width=.9\linewidth]{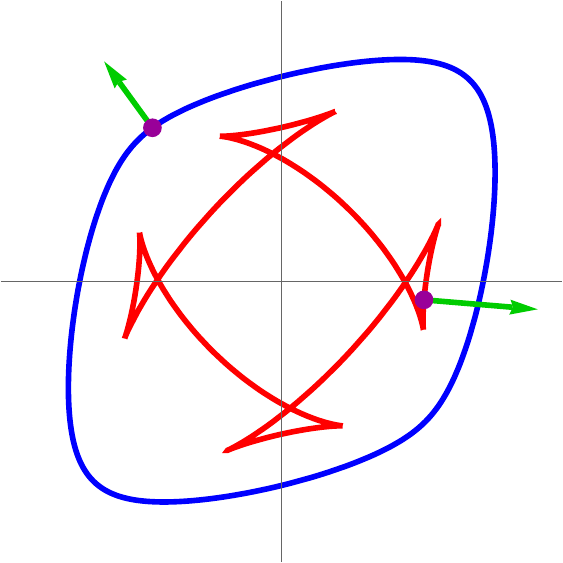}
  \caption{The set of group velocities in $T_x\R^2$.
  The group velocity corresponding to $\vect{p}\in T^*_x\R^2$ with $\lambda_i(\Gamma(\vect{p}))=1$ is the gradient $\nabla\lambda_i(\Gamma(\vect{p}))\in T_x\R^2$.
  The group velocity describes the actual movement of the phonons.
  The phase velocity is $\vect{p}/\abs{\vect{p}}^2$ defined with the Euclidean norm.
  }
  \label{fig:sub2}
\end{subfigure}
\caption{A 2-dimensional slowness curve $S_x$ with $x\in\R^2$ fixed and parameters
$b_{11}=10$,
$b_{22}=12$,
$b_{33}=20$,
$b_{12}=2$,
$b_{23}=5$,
$b_{13}=3$; cf.~\S\ref{sec:alg-principles}.
The slowness curve is on the left and the group velocity curve is on the right.
Here, the multiplicity of each eigenvalue in the Christoffel matrix is one for all $\vect{p}\neq0$, so
different eigenvectors (polarizations) correspond to different eigenvalues~$\lambda_i$.
The slowness curve $S_x$ is the set where at least one of the~$\lambda_i$'s equals~$1$.
The bigger eigenvalue corresponds to 
\textcolor{blue}{quasi-pressure polarization (qP)} and the smaller one to \textcolor{red}{quasi-shear (qS)}.
The \textcolor{color1}{gradients of the eigenvalue functions on the left} correspond to \textcolor{color1}{points on the right} and \textcolor{color2}{vice versa}, via Legendre duality.
The \textcolor{blue}{blue curve} on the right is the unit sphere of the Finsler geometry whose geodesics the wave packets follow. Polarization vectors are not indicated in these figures.
For a description of how the slowness curve arises in geophysics and microlocal analysis see, e.g., Yedling~\cite{Yedling}, who provides a geometric analysis of the wave front in a homogeneous anisotropic material.
}
\label{fig:2d-slowness-example}
\end{figure}

The set of points
\begin{equation*}
S_x = \left\{ \vect{p}\in T^*_x\R^n : \det( \Gamma(x,\vect{p})-I_n ) = 0 \right\}
\end{equation*}
is called the \defi{slowness (hyper)surface} at the point~$x$. A point $\vect{p}\in T^*_x\R^n$ belongs to the slowness surface exactly when~$1$ is an eigenvalue of $\Gamma(x,\vect{p})$; since the Christoffel matrix is $2$-homogeneous in~$\vect{p}$, the slowness surface encodes the eigenvalue information of the Christoffel matrix.
See Figure~\ref{fig:2d-slowness-example} for an example of a slowness surface when $n=2$, where $S_x$ is a curve.

\emph{A priori}, the slowness surface contains no information about the eigen\emph{vectors} of the Christoffel matrix. These vectors are the polarizations of singularities and correspond to the direction of oscillation, whereas the slowness vector corresponds roughly to the direction of propagation. Our method is suited for situations where we observe only the singularities in space but not their polarizations. Singularities without polarization can aptly be called \defi{unpolarized phonons}; the phonon is the particle (or wave packet) corresponding to the displacement field in wave--particle duality; cf.~\cite{W:phonon-book}. The eigenvalues of the Christoffel matrix give rise to Hamiltonians --- one for each polarization --- that determine the time evolution of unpolarized phonons. The slowness surface is the union of the unit level sets of these Hamiltonians; it might not split cleanly into~$n$ branches and~$n$ well-defined and smooth Hamiltonians, due to degenerate eigenvalues, so we treat the slowness surface as a single object.

\subsubsection{An algebraic view}
\label{ss:algebraicview}

Considering the \defi{slowness polynomial} $\spoly{\vect{a}}(\vect{p})\coloneqq\det(\Gamma(x,\vect{p})-I_n)$ as a polynomial of degree $2n$ in the variables $p_1,\dots,p_n$, and parametrized by $\vect{a}\in\stiffnesstensor{n}$, the slowness surface is an object of a very algebraic nature. Our core algebraic result, Theorem~\ref{thm:GenericIrreducibility} is that for generic~$\vect{a}$, the slowness surface is an \emph{irreducible algebraic variety}. The map that takes a stiffness tensor into the corresponding slowness polynomial,
\begin{equation}
\label{eq:Fmap}
\begin{split}
\spolymap \colon\stiffnesstensor{n}&\to\R[p_1,\dots,p_n] ,\\
\vect{a} &\mapsto \spoly{\vect{a}},
\end{split}
\end{equation}
can be used to formulate many of our key statements.

As our focus is on the analysis on a fiber of the cotangent bundle rather than the whole bundle, from now on we consider the point $x\in\R^n$ fixed and drop it from the notation. Thus, for example, the Christoffel matrix will henceforth be denoted $\Gamma(\vect{p})$. Similarly, we call the reduced (or density-normalized) stiffness tensor $\vect{a}=\rho^{-1}\vect{c}$ the stiffness tensor for simplicity.

%%%%%%%%%%%
\subsection{Departure point}

The goal of a typical inverse problem in anisotropic linear elasticity is to reconstruct the stiffness tensor field from some kind of boundary measurement --- or to prove that the field is uniquely determined by ideal boundary data. A good first step is to analyze the propagation of singularities microlocally from travel time data or hyperbolic Cauchy data\footnote{The hyperbolic Cauchy data is the set of all Dirichlet and Neumann boundary values of all solutions to the elastic wave equation. It is the graph of the Dirichlet-to-Neumann map.}. This turns the analytic inverse problem into a geometric one, where the task is to recover the geometry governing the propagation of singularities.

In the anisotropic elastic setting this geometry is quite complicated. The innermost branch of the slowness surface, called the qP or \defi{quasi-pressure branch} (see Figure~\ref{fig:2d-slowness-example}), determines a Finsler geometry when the highest eigenvalue of the Christoffel matrix is non-degenerate; see e.g.~\cite{dHILS:BDF} for details. Finding the Finsler function on some subset of the tangent bundle amounts to finding a subset of the slowness surface at some or all points. For other polarizations (qS), the unit cosphere --- which is a branch of the slowness surface --- may fail to be convex. It is also common for the slowness surface to have singular points where two branches meet (see \cite{I:degenerate} for details), which is an obstruction to having a smooth and globally defined Finsler geometry for slower polarizations.

Despite these issues, at most points $x \in \R^n$ and in most directions $\vect{p} \in T_x^*\R^n$ the Christoffel matrix $\Gamma(x,\vect{p})$ has~$n$ different eigenvalues. In the neighborhood of such a point $(x,\vect{p})$ on the cotangent bundle the elastic wave equation~\eqref{eq:EWE} can be diagonalized microlocally and any solution splits nicely into~$n$ different polarizations. It is also this non-degenerate setting where our description of propagation of singularities is valid without additional caveats.

Ideally, the solution of such a geometric inverse problem produces a qP Finsler geometry in full. The unit cosphere of this geometry at each point is the qP branch of the slowness surface, e.g., \cite{dHILS:BDF,dHILS:BSR}. In some cases the recovery is not full but only a part of the slowness surface can be reconstructed; see \cite{dHILS:BDF,dHILS:BSR,dHIL:Finsler-Dix}.  In several applications one measures only the arrival times of the fastest waves which give the travel times related to the qP polarized waves. Further investigation will surely lead to mathematical results that provide full or partial information of the slowness surface of a single polarization. For our purposes, it is irrelevant where the partial knowledge of the slowness surface comes from, only that this information is indeed accessible.

Our contribution is to take the next step: We prove that generically a small subset of one branch\footnote{It is unimportant whether the small patch of the slowness surface we start with is qP or some other polarization. We mainly refer to qP only because it is the easiest polarization to measure physically and mathematically, corresponding to the fastest waves and a well-behaved Finsler geometry.} of the slowness surface determines uniquely the entire slowness surface (with all branches) and, if possible, the stiffness tensor field. No polarization information is required as input for our methods. Taking the determinant of the principal symbol amounts to ignoring polarization information. However, we can fill in the polarizations of the singularities \emph{after} reconstructing the stiffness tensor field.

%%%%%%%%%%%
\subsection{Algebraic goals}

The~$n$ different branches of the slowness surface, each corresponding at least locally to a polarization, are coupled together by the simple fact that they are all in the vanishing set of the same polynomial --- the slowness polynomial. If we can show that this polynomial is irreducible (i.e., cannot be written as a product of two polynomials in a non-trivial way), then a small open subset of the slowness surface can be completed to the whole slowness surface by taking its Zariski closure, i.e., taking the vanishing set of a polynomial of the right shape that interpolates the known small set. Physically, this means that a small (Euclidean) open subset of the slowness surface, e.g., an open subset of the sheet of the slowness surface associated to qP-waves, determines both the entire sheet, as well as the sheets of the slowness surface associated to qS-polarized waves. Even though the sheet associated to qP-polarized waves and the sheets associated to qS-polarized waves
are disjoint in the Euclidean topology, we will show that for stiffness tensors in a generic set these sheets all lie in the same connected component in the complex Zariski topology of the slowness surface.

We note that determining all the sheets of a slowness surface, or the stiffness tensor that gave rise to surface, from a small open subset of the qP-sheet, is impossible for \emph{some} stiffness tensors. For example, fully isotropic stiffness tensors parametrized by the two Lam\'e parameters, which describe many real materials, give rise to the slowness polynomials of the form
\begin{equation}
\label{eq:slownesspolyfullyisotropic}
P(\vect{p}) = (c_P^{2}\abs{\vect{p}}^2-1) (c_S^{2}\abs{\vect{p}}^2-1)^{n-1},
\end{equation}
where $c_P$ and~$c_S$ are the pressure and shear wave speeds of the material respectively. Such a polynomial is manifestly reducible. The slowness surface for such materials consists of two concentric spheres, the inner one of which is degenerate if $n>2$. The two radii are independent. If we know a small subset of the pressure sphere, taking the Zariski closure completes it into the whole sphere but contains no information about the other sphere. This, however, is an exceptional property of a highly symmetric stiffness tensor. We therefore set out to prove that the slowness polynomial is irreducible for \emph{most} stiffness tensors.

Once the whole slowness surface or slowness polynomial is recovered, we can use it in practice to reconstruct the stiffness tensor. The reason why this works is less obvious, but it will be explained in~\S\ref{sec:alg-principles} once we lay out some preliminaries.

From the point of view of inverse problems in analysis, geometry, and linear elasticity, we need two algebraic results that are straightforward to state but less straightforward to prove. These key results are given in~\S\ref{sec:main-results} below, with the definitions for Theorem~\ref{thm:2-layers} detailed in~\S\ref{sec:2-layer-def}. The algebraic results mentioned below hinge on the more technical results given in~\S\ref{sec:results}.

%%%%%%%%%%%
\subsection{Main results}
\label{sec:main-results}

We have two main results on inverse problems. The first one states that the whole slowness surface, including all the branches, is uniquely determined by any small open set of a single branch. That is, partial data for a single polarization is enough for full uniqueness.

\begin{mainthm}[Uniqueness of slowness surfaces and stiffness tensor from partial data]
\label{thm:generic-stiffness-uniqueness}
\ 
\begin{enumerate}[leftmargin=*]
    \item \label{item:thmA-surface}
    Let the dimension of the space be $n\in\{2,3\}$.
    There is an open and dense subset $U\subset\posstiffnesstensor{n}$ of the set of positive stiffness tensors such that the following holds:
    If $\vect{a}\in U$, then any non-empty Euclidean relatively open subset of the slowness surface corresponding to~$\vect{a}$ determines the whole slowness surface uniquely.
    \item \label{item:thmA-tensor}
    Let the dimension of the space be $n\in\{2,3\}$ and if $n=3$ suppose that Conjecture~\ref{conj:uniquereconstruction} holds.
    There is an open and dense subset $U\subset\posstiffnesstensor{n}$ of the set of positive stiffness tensors such that the following holds:
    If $\vect{a}\in U$, then any non-empty Euclidean relatively open subset of the slowness surface corresponding to~$\vect{a}$ determines the stiffness tensor~$\vect{a}$ uniquely.
\end{enumerate}
\end{mainthm}

The conclusion of Theorem~\ref{thm:generic-stiffness-uniqueness} is false for isotropic stiffness tensors in $E^+(n)$, because the two wave speeds (radii of spheres) in~\eqref{eq:slownesspolyfullyisotropic} are independent. Theorem~\ref{thm:generic-stiffness-uniqueness} makes precise the idea that this failure is a rare exception among all stiffness tensors. The openness and density of the sets of stiffness tensors in Theorem~\ref{thm:generic-stiffness-uniqueness} is valid both in the Euclidean topology and the Zariski topology.

The next theorem states roughly that, generically, a two-layer model of a planet with piecewise constant stiffness tensor field is uniquely determined by geometric travel time data for rays traversing the interior of the planet. The two layers are a highly simplified model of the Earth with a mantle and a core, both with homogeneous but anisotropic materials.

We consider two data types for a two-layer model where in the outer layer $\Omega\setminus \omega\subset \R^n$ the stiffness tensor is equal to $\vect{A}$
and in the inner layer $\omega\subset \R^n$ the stiffness tensor is equal to $\vect{a}$. The first data type, denoted $\nakeddata=\nakeddata(\omega,\Omega,\vect{a},\vect{A})$, contains only travel time information between boundary points, while the second data type, denoted $\decorateddata=\decorateddata(\omega,\Omega,\vect{a},\vect{A})$, contains also directional information at the boundary. The precise definitions of these data sets are given below in~\S\ref{sec:2-layer-def}.
By $X_1\approx X_2$ for two subsets~$X_i$ of the same Euclidean space, we mean that there are dense open subsets $U_i\subset X_i$ for $i=1,2$ so that $U_1=U_2$. See~\S\ref{sec:2-layer-def} for details of the various kinds of rays and data, and the notion of admissibility. 

\begin{mainthm}[Two-layer model]
\label{thm:2-layers}
Let $n\in\{2,3\}$ and if $n=3$ suppose that Conjecture~\ref{conj:uniquereconstruction} holds.
For both $i=1,2$, let $\omega_i$, $\Omega_i\subset\R^n$ be open domains satisfying
$\overline \omega_i \subset\Omega_i$ and
$\Omega_1=\Omega_2\eqqcolon\Omega$.
There is an open and dense subset~$U$ in the space of admissible pairs of stiffness tensors so that for all $(\vect{a}_i,\vect{A}_i)\in U$ the following holds.

\smallskip

If
$\nakeddata(\omega_1,\Omega,\vect{a}_1,\vect{A}_1)
\approx
\nakeddata(\omega_2,\Omega,\vect{a}_2,\vect{A}_2)$,
then $\vect{A}_1=\vect{A}_2$ and $\omega_1=\omega_2$.

\smallskip

If
$\decorateddata(\omega_1,\Omega,\vect{a}_1,\vect{A}_1)
\approx
\decorateddata(\omega_2,\Omega,\vect{a}_2,\vect{A}_2)$,
then $\vect{A}_1=\vect{A}_2$, $\omega_1=\omega_2$, and additionally $\vect{a}_1=\vect{a}_2$.
\end{mainthm}

\begin{remarks}\ 
\begin{enumerate}[leftmargin=*]
\item The conclusions of Theorem~\ref{thm:2-layers} are valid only for generic stiffness tensors. They are false, for example, for all isotropic stiffness tensors.
\medskip
    \item The equivalence $\approx$ ensures that exceptional rays play no role --- possible exceptions include gracing rays, zero transmission or reflection coefficients, and cancellation after multipathing. All these issues are typically rare; e.g. the transmission and reflection coefficients are in many cases analytic functions with isolated zeros~\cite{Heijden}. If one assumes that the full data sets are equal, then one might prove results like ours by only comparing which rays are missing from the data or behave exceptionally. To ensure that the conclusion is reached using well-behaved rays and that the omission or inclusion of a small set of rays (whether well or ill behaved) is irrelevant, we will only assume that the data sets are only almost equal.
\medskip
\item The second part of the statement should be seen in light of the first one. Namely, the second assumption implies the first one, so the stiffness tensor~$\vect{A}$ is uniquely determined by the data. Given this stiffness tensor in the outer layer, the knowledge of the slowness vector is almost the same as knowing the velocity vector or the tangential component or length of either one. Full directional data can be used to find which slowness vectors are admissible, and that in turn generically determines the stiffness tensor. The exact details are unpleasant, so the statement of Theorem~\ref{thm:2-layers} has been optimized for readability rather than strength.
The result as stated above can be adapted to other measurement scenarios.
\medskip
\item All polarizations are included in the data of Theorem~\ref{thm:2-layers}.
The proof is based only on the fastest one (qP), but ignoring the other ones is not trivial. Even if incoming and outgoing waves at the surface are qP, there can be segments in other polarizations due to mode conversions inside.
\medskip
\item Theorem~\ref{thm:2-layers} was stated geometrically.
Geometric data of this kind may be obtained from boundary data for the elastic wave equation~\eqref{eq:EWE}. It is not uncommon in geophysics to work directly with geometric ray data; see e.g.~\cite{raytheory-example}.
\end{enumerate}
\end{remarks}

\begin{remark}
Theorems~\ref{thm:generic-stiffness-uniqueness} and~\ref{thm:2-layers} above are conditional on Conjecture~\ref{conj:uniquereconstruction} in 3D for the determination of the stiffness tensor. If we know what that the stiffness tensor is monoclinic or orthorhombic and we know the symmetry planes, then we get a similar result unconditionally by Theorem~\ref{thm:polynomial-to-tensor}. Uniqueness in these symmetry classes, however, is local instead of global due to the nature of Theorem~\ref{thm:polynomial-to-tensor}: every stiffness tensor in the generic set has a small neighborhood (in the Euclidean topology) where there are no anomalous companions (other stiffness tensors with the same slowness surface). If the symmetry planes are not known, they first have to be determined from the slowness surface itself.
\end{remark}

These results on inverse problems hinge on two algebraic results. An analytically oriented reader should think that the slowness polynomial contains the same information as the corresponding slowness surface, and this is indeed the case, set-theoretically.

\begin{mainthm}[Generic irreducibility]
\label{thm:GenericIrreducibility}
The slowness polynomial associated to a generic stiffness tensor in dimension $n \in \{2,3\}$ is irreducible over~$\C$. The same conclusion holds for slowness polynomials arising from generic orthorhombic and generic monoclinic stiffness tensors when $n = 3$.
\end{mainthm}

The word \defi{generic} in Theorem~\ref{thm:GenericIrreducibility} is used in the sense of algebraic geometry: Let~$m$ be the number of distinct components of a reduced stiffness tensor. A generic set of stiffness tensors is a subset of $\R^m$ whose complement is a finite union of algebraic subsets of~$\R^m$ of dimension $\leq m-1$, each of which is defined by a finite set of polynomials. In our case, this complement parametrizes the collection of stiffness tensors giving rise to reducible slowness surfaces; it is not empty: we already know that the slowness polynomial~\eqref{eq:slownesspolyfullyisotropic} associated to a fully isotropic stiffness tensor is reducible. A priori, it is possible that the sets $\orthorhombic$ and $\monoclinic$ belong to this complement. Thus, the second sentence in Theorem~\ref{thm:GenericIrreducibility} is not simply a consequence of the first.

In the spirit of modern algebraic geometry, we prove Theorem~\ref{thm:GenericIrreducibility} by considering \emph{all} slowness polynomials at once, in a family
\[
f \colon \bS \to \Aff^k_\R,
\]
where the coordinates at a point $y = y(x) \in \R^k = \Aff^k(\R)$ record the coefficients of a single slowness polynomial at $x \in \R^n$, and the corresponding fiber $f^{-1}(y) \subset \bS$ is the slowness surface~$S_x$. The principle of \emph{generic geometric integrality}, due to Grothendieck, ensures that if the map~$f$ satisfies a few technical hypotheses, then there is a Zariski open subset of $Y \subset \Aff_\R^k$ such that the individual slowness surfaces in $f^{-1}(Y)$ are irreducible, even over~$\C$ (equivalently, their corresponding slowness polynomials are $\C$-irreducible). A Zariski open subset of~$\Aff_\R^k$ is dense for the Euclidean topology, as long as it is not empty.  Thus, we must check by hand the existence of a single $\C$-irreducible slowness polynomial to conclude. In the case $n=3$, we use an explicit stiffness tensor, modelling a specific physical mineral, olivine, to verify the non-emptiness of the set~$Y$. This task is accomplished by reduction to modulo a suitably chosen prime (see Lemma~\ref{lem:IrreducibilityCriterion}---we include a proof for lack of a reference to this tailored lemma; we hope it will be useful in other inverse-problem contexts). 

Our second algebraic result explores the correspondence~\eqref{eq:Fmap} between stiffness tensors and slowness polynomials. We show that this correspondence is generically injective in dimension two, and we conjecture the same holds in dimension three (Conjecture~\ref{conj:uniquereconstruction}). In contrast, Helbig and Carcione in~\cite{HC2009} have already observed that a slowness polynomial associated to an orthorhombic stiffness tensor can arise from more than one such tensor. They gave sufficient conditions for the existence of what they called \defi{anomalous companions} of an orthorhombic stiffness tensor.  We tighten their results to show that their conditions are also generically necessary. In analogy with their terminology, we show that a generic monoclinic stiffness tensor has precisely one anomalous companion that gives rise to the same slowness polynomial.

\begin{mainthm}[Generic Reconstruction]
\label{thm:polynomial-to-tensor}
\ 
\begin{enumerate}[leftmargin=*]
    \item \label{item:2D}
    The map $\spolymap\colon \stiffnesstensor{2} \to \R[p_1,p_2]$ is generically injective; i.e., in dimension $n = 2$, the slowness polynomial associated to a generic stiffness tensor determines the stiffness tensor \emph{uniquely}.
    \smallskip
    \item \label{item:3Dortho} The map $\spolymap|_{\orthorhombic} \colon \orthorhombic \to \R[p_1,p_2,p_3]$ is generically $4$-to-$1$; i.e., a generic orthorhombic stiffness tensor is determined by its slowness polynomial up to three anomalous companions.
    \smallskip
    \item \label{item:3Dmono}
    The map $\spolymap|_{\monoclinic} \colon \monoclinic \to \R[p_1,p_2,p_3]$ is generically $2$-to-$1$; i.e., a generic monoclinic stiffness tensor is determined by its slowness polynomial up to one anomalous companion.
\end{enumerate}
\end{mainthm}

We give \emph{three} proofs of Theorem~\ref{thm:polynomial-to-tensor}(\ref{item:2D}), showcasing multiple algebro-geometric approaches to this inverse problem.
The first proof is a direct application of a result of Eisenbud and Ulrich~\cite{EisenbudUlrich}. The second proof is more ``hands-on'' and generalizes well to parts~(\ref{item:3Dortho}) and~(\ref{item:3Dmono}) of the theorem. Finally, the third proof ensues from studying the following question: Given a polynomial with real coefficients, what conditions must its coefficients satisfy for it to be a slowness polynomial? In other words: can we characterize slowness polynomials among all polynomials?  The methods we employ when $n=2$ are in principle applicable when $n=3$ as well. The required computer calculations, however, are currently infeasible, causing us to restrict to specific symmetry types of stiffness tensors. With enough computing power, our methods will produce a proof of Theorem~\ref{thm:polynomial-to-tensor}(\ref{item:2D}) when $n=3$ --- if it is true.

\begin{mainconj}
\label{conj:uniquereconstruction}
The slowness polynomial associated to a generic stiffness tensor in dimension $n = 3$ determines the stiffness tensor; i.e., the map $\spolymap \colon \stiffnesstensor{3} \to \R[p_1,p_2,p_3]$ is generically injective.
\end{mainconj}

Calculations on specific slowness polynomials support Conjecture~\ref{conj:uniquereconstruction}: for any given slowness polynomial, we can use Gr\"obner bases to instantaneously reconstruct the set of all possible stiffness tensors associated with the polynomial (see~\S\ref{ss:2DGB}). For a concrete example of a fully unique stiffness tensor, we carry out this unique reconstruction process for the slowness polynomial of albite~\cite{albite}, an abundant mineral in the Earth's crust; see~\cite{this-paper}.

Theorems~\ref{thm:GenericIrreducibility} and~\ref{thm:polynomial-to-tensor} are proved in~\S\ref{sec:AG-proofs}; however, for the benefit of readers without much exposure to algebraic geometry, we explain the principles involved in~\S\ref{sec:alg-principles} in the case $n = 2$, to avoid clutter.

\subsubsection{Positivity}

A physical stiffness tensor has to be positive, but its anomalous companions may fail to be so. In the orthorhombic setup the anomalous companions and their positivity can be neatly described via the Cayley cubic surface, a central object in the classical algebro-geometric canon; see \S\ref{sec:Cayley-cubic} for details. It would be quite interesting to perform a similar positivity analysis for the anomalous companion of a monoclinic tensor.

%%%%%%%%%%%
\subsection{Related results}

Motivated by seismological considerations, inverse boundary value problems in elasticity have been studied since 1907, when Wiechert and Zoeppritz posed them in  their paper ``\"Uber Erdbebenwellen''(On Earthquake Waves) \cite{Wiechert}; see also \cite{Ammari,Beretta,PSU}. The first breakthrough results in \textit{elastostatics} for isotropic media were by Nakamura and Uhlmann~\cite{NakamuraUhlmann_1}, followed by results by Eskin and Ralston~\cite{EskinRalston_1} for full boundary data and Imanuvilov, Uhlmann and Yamamoto~\cite{IUY_1} for partial boundary data. Stefanov, Uhlmann and Vasy~\cite{SUV} studied recovery of smooth \textit{P}- and \textit{S}-wave speeds in the elastic wave equation from knowledge of the Dirichlet-to-Neumann map in the isotropic case, see also  \cite{Barcelo} on the reconstruction of the density tensor. Beretta, Francini and Vessella~\cite{Beretta2} studied the stability of solutions to inverse problems. Uniqueness results for the tomography problem with interfaces, again, in the isotropic case, in the spirit of Theorem~\ref{thm:2-layers}, were considered by Caday, de Hoop, Katsnelson and Uhlmann~\cite{CdHKU}, as well as by Stefanov, Ulhmann and Vasy~\cite{SUVII}. 

The related inverse travel time problem (for the corresponding Riemannian metric) has been studied in isotropic media using integral geometry in \cite{Michel,SU,SUV1,SUV2,UV} and metric geometry in \cite{Burago}.

Anisotropic versions of the \textit{dynamic} inverse boundary value problem have been studied in various different settings. Rachele and Mazzucato studied the geometric invariance of elastic inverse problems in  \cite{Mazzucato1}.  In \cite{RacheleMazzucato_2007,Mazzucato2}, they showed that for certain classes of transversely isotropic media, the slowness surfaces of which are ellipsoidal, two of the five material parameters are partially determined by the dynamic Dirichlet-to-Neumann map. Before that, Sacks and Yakhno \cite{SacksYakhno_1998} studied the inverse problem for a layered anisotropic half space using the Neumann-to-Dirichlet map as input data, observing that only a subset of the components of the stiffness tensor can be determined expressed by a ``structure'' condition.
De Hoop, Nakamura and Zhai \cite{dHNakamuraZhai_2019} studied the recovery of piecewise analytic density and stiffness tensor of a three-dimensional domain from the local dynamic Dirichlet-to-Neumann map. They give global uniqueness results if the material is transversely isotropic with known axis of symmetry or orthorhombic with known symmetry planes on each subdomain. They also obtain uniqueness of a fully anisotropic stiffness tensor, assuming that it is piecewise constant and that the interfaces which separate the subdomains have curved portions. Their method of proof requires the use of the (finite in time) Laplace transform.
Following this transform, some of the techniques are rooted in the proofs of analogous results for the inverse boundary value problem in the elastostatic case \cite{NakamuraTUhlmann_1999, CHondaNakamura_2018}. C\^{a}rstea, Nakamura and Oksanen \cite{CNakamuraOksanen_2020} avoid the use of the Laplace transform and obtain uniqueness, in the piecewise constant case, closer to the part of the boundary where the measurements are taken for shorter observation times and further away from that part of the boundary for longer times.

Under certain conditions, the dynamic Dirichlet-to-Neumann map determines the scattering relation, allowing a transition from analytic to geometric data. Geometric inverse problems in anisotropic elasticity have received increasing attention over the past few years. In the case of transversely anisotropic media the elastic parameters are determined by the boundary travel times of all the polarizations \cite{dHUhlmannVasy_2020,Zou_2021}. A compact Finsler manifold is determined by its boundary distance map~\cite{dHILS:BDF}, a foliated and reversible Finsler manifold by its broken scattering relation~\cite{dHILS:BSR}, and one can reconstruct the Finsler geometry along a geodesic from sphere data~\cite{dHIL:Finsler-Dix}.

Linearizing about the isotropic case, that is, assuming ``weak'' anisotropy, leads to the mixed ray transform for travel times between boundary points. De Hoop, Saksala, Uhlmann and Zhai \cite{dHUSZ:MRT} proved ``generic'' uniqueness and stability for this transform on a three-dimensional compact simple Riemannian manifold with boundary, characterizing its kernel. Before that De Hoop, Saksala and Zhai \cite{dHSaksalaZhai_2019} studied the mixed ray transform on simple $2$-dimensional Riemannian manifolds. Linearizing about an isotropic case but only with conformal perturbations leads to a scalar geodesic ray transform problem on a reversible Finsler manifold, and the injectivity of that transform was established in~\cite{IM:Finsler-XRT} in spherical symmetry.

Assuming lack of symmetry naturally leads to the occurrence of singular points in the slowness surface. This is inherent in exploiting algebraic geometry to obtain the results in this paper. However, the singular points lead to fundamental complications in the application of microlocal analysis to a parametric construction revealing the geometry of elastic wave propagation, see \cite{Greenleaf2,Greenleaf3,Greenleaf1}. The points are typically associated with conical refraction \cite{MelroseUhlmann_1979, Uhlmann_1982, Dencker_1988, BraamDuistermaat_1993, ColindeVerdiere_2003, ColindeVerdiere_2004}.

%%%%%%%%%%%
\subsection{Outline}

The paper is organized as follows. In~\S\ref{sec:alg-principles} we explain the general algebro-geometric framework underlying the proofs of Theorems~\ref{thm:GenericIrreducibility} and~\ref{thm:polynomial-to-tensor} in dimension~$2$, where the number of parameters is small, making it easier to digest the ideas involved. We pivot in \S\S\ref{sec:2-layer-def}--\ref{sec:IP-proofs} to the study of inverse problems, setting up precise definitions for Theorems~\ref{thm:generic-stiffness-uniqueness} and~\ref{thm:2-layers} in~\S\ref{sec:2-layer-def} and giving proofs for these theorems in~\S\ref{sec:IP-proofs}. In~\S\ref{sec:AG-proofs} we prove Theorems~\ref{thm:GenericIrreducibility} and~\ref{thm:polynomial-to-tensor}.

%%%%%%%%%%%
\subsection*{Acknowledgements}
We thank Mohamed Barakat, Olivier Benoist, Daniel Erman, Bjorn Poonen, and Karen Smith for useful discussions around the algebro-geometric content of the paper.  We are particularly indebted to Jan Draisma, who pointed out an error in an earlier version of this paper.
M.V.\ de~H.\ was supported by the Simons Foundation under the MATH + X program, the National Science Foundation under grant DMS-2108175, and the corporate members of the Geo-Mathematical Imaging Group at Rice University.
J.\ I.\ was supported by the Research Council of Finland (Flagship of Advanced Mathematics for Sensing Imaging and Modelling grant 359208; Centre of Excellence of Inverse Modelling and Imaging grant 353092; and other grants 351665, 351656, 358047) and the V\"ais\"al\"a project grant by the Finnish Academy of Science and Letters.
M.\ L.\ was supported by the European Research Council of the European Union, the ERC-AdG project PDE-Inverse, 101097198, and Academy of Finland grants 284715 and 303754, the Centre of Excellence of Inverse Modelling and Imaging and the Flagship of Advanced Mathematics for Sensing Imaging and Modelling.
A.\ V.-A.\ was partially supported by NSF grants DMS-1902274 and DMS-2302231,  as well as NSF Grant DMS-1928930 while in residence at MSRI/SLMath in Berkeley (Spring 2023). He thanks Fran\c cois Charles for hosting him at the D\'epartement de Math\'ematiques et Applications of the \'Ecole Normale Sup\'erieur in Summer 2022, where part of this paper was written.
Views and opinions expressed are those of the authors only and do not necessarily reflect those of the European Union or the other funding organizations.

%%%%%%%%%%%%%%%%%%%%%%%%%%%%%%%%
\section{Algebro-geometric principles: a case study in dimension 2}
\label{sec:alg-principles}

To illustrate how algebraic geometry bears on inverse problems in anisotropic elasticity, we consider a two-dimensional model, where a slowness surface is in fact a curve.  An anisotropic stiffness tensor in this case is determined by six general real parameters: although such a tensor $\vect{a} = (a_{ijkl})\in \R^{16}$ has $16$ components, once we take into account the major and minor symmetries of the tensor, only $6$ distinct parameters are left.
Following Voigt's notation~(see~\S\ref{sec:tensor-alg}), they are
\begin{align*}
    b_{11} &= a_{1111}, & b_{22} &= a_{2222}, \\
    b_{12} &= a_{1122} = a_{2211}, & b_{23} &= a_{2212} = a_{2221} = a_{1222} = a_{2122}, \\
    b_{13} &= a_{1112} = a_{1121} = a_{1211} = a_{2111}, & b_{33} &= a_{1212} = a_{2112} = a_{1221} = a_{2121}.
\end{align*}
The corresponding Christoffel matrix is
\begin{equation}
\label{eq:2DChristoffelmatrix}
\Gamma(\vect{p}) =
\begin{pmatrix}
b_{11} p_1^2 + 2b_{13} p_1p_2 + b_{33} p_2^2 & b_{13}p_1^2 + (b_{33}+b_{12}) p_1p_2 + b_{23} p_2^2 \\
b_{13}p_1^2 + (b_{33}+b_{12}) p_1p_2 + b_{23}p_2^2 & b_{33}p_1^2 + 2b_{23}p_1p_2 + b_{22}p_2^2
\end{pmatrix}.
\end{equation}
The \defi{slowness curve} is the vanishing set of $P_{\vect{a}}(\vect{p}) := \det\left(\Gamma(\vect{p}) - I_2\right)$ in $\R^2$:
\[
S \coloneqq \{\vect{p} \in \R^2 \mid P_{\vect{a}}(\vect{p}) = 0\}.
\]
The polynomial $P_{\vect{a}}(\vect{p})$ has degree~$4$ in the variables $p_1$ and~$p_2$, but not every monomial of degree $\leq 4$ in $p_1$ and~$p_2$ appears in it.  In fact, 
\begin{equation}
\label{eq:gen2Dpoly}
P_{\vect{a}}(\vect{p}) = c_1p_1^4 + c_2p_1^3p_2 + c_3p_1^2p_2^2 + c_4p_1^2 + c_5p_1p_2^3 + c_6p_1p_2 + c_7p_2^4 + c_8p_2^2 + c_9,
\end{equation}
for some constants $c_i \in \R$, $i = 1\dots,9$.  These constants are not arbitrary; they have to satisfy relations like $c_9 = 1$, or the more vexing
\begin{equation}
\label{eq:oneconditiononcs}
\begin{split}
-4c_1^2 &+ 4c_1c_4c_8 - c_1c_6^2 + 8c_1c_7 - 4c_1c_8^2 - c_2^2 - 2c_2c_5 \\
&\quad + 2c_2c_6c_8 - c_3c_6^2 - 4c_4^2c_7 + 2c_4c_5c_6 + 4c_4c_7c_8 - c_5^2 - c_6^2c_7 - 4c_7^2 = 0 
\end{split}
\end{equation}
in order to arise from a stiffness tensor (see \S\ref{sec:uniquereconstruction} and \S\ref{sec:slownesspolys}).

%%%%%%%%%%%
\subsection{Goals}
We want to know two things: 
\begin{enumerate}
    \item For generic choices of the parameters~$b_{ij}$, the curve $S \subset \R^2$ is irreducible, even over the complex numbers.
    \medskip
    \item For generic choices of $c_1,\dots,c_8$ corresponding to a slowness polynomial, there is a unique set of $b_{ij}$'s giving rise to the polynomial~\eqref{eq:gen2Dpoly}, and this polynomial can be explicitly computed if we approximate $c_1\dots,c_8$ by rational numbers. 
\end{enumerate}
We accomplish both of these goals by leveraging powerful results in both the theory of schemes\footnote{Schemes over a field form a category that is richer and more flexible than the corresponding category of varieties.}, as developed by Alexander Grothendieck, and the application of computational techniques under the banner of Gr\"obner bases.

%%%%%%%%%%%
\subsection{Generic Irreducibility}
To realize our first goal, we must compactify a slowness curve and consider \emph{all slowness curves at once}, in a family.  This allows us to apply a suite of results from scheme theory, including ``general geometric integrality''. So think now of the parameters~$b_{ij}$ as indeterminates, and the entries of $\Gamma(\vect{p})$ as belonging to the polynomial ring $A[p_0,p_1,p_2]$ with coefficients in the polynomial ring $A \coloneqq \R[b_{11},b_{12},b_{13},b_{22},b_{23},b_{33}]$. Then the \defi{homogenized slowness polynomial}\footnote{Perhaps $\tilde P_{\vect{a}}(\vect{p})$ would be more appropriate notation than $\tilde P(\vect{p})$, but we drop the $\vect{a}$ to avoid clutter.} is given by
\begin{equation}
\label{eq:slownesspoly2D}
\begin{split}
\tilde P(\vect{p}) &= \det(\Gamma(\vect{p})-p_0^2I_2) \\
&= (b_{11}b_{33}-b_{13}^2)p_1^4 + 2(b_{11}b_{23}-b_{12}b_{13})p_1^3p_2 
    + (b_{11}b_{22}-b_{12}^2-2b_{12}b_{33}+2b_{13}b_{23})p_1^2p_2^2 \\
&\quad - (b_{11}+b_{33})p_1^2p_0^2 + 2(b_{13}b_{22}-b_{12}b_{23})p_1p_2^3 
    - 2(b_{13}+b_{23})p_1p_2p_0^2 + (b_{22}b_{33}-b_{23}^2)p_2^4 \\
&\quad - (b_{22}+b_{33})p_2^2p_0^2 + p_0^4.
\end{split}
\end{equation}
It is a homogeneous polynomial of degree~$4$ in the variables $p_0$, $p_1$ and $p_2$, with coefficients in the polynomial ring $A$. Its zero-locus traces a curve in the projective plane $\PP^2_A$ with homogeneous coordinates $(p_0 : p_1 : p_2)$ and coefficient ring $A$. Since $\PP^2_A$ is naturally isomorphic, as an $\R$-scheme, to the fibered product $\Aff^6_\R \times \PP^2_\R$, the vanishing of $\tilde P(\vect{p})$ is also naturally a hypersurface in the product of the affine space over $\R$ with coordinates $b_{11},\dots,b_{33}$ and the projective plane over $\R$ with homogeneous coordinates $(p_0 : p_1 : p_2)$. Let
\begin{align*}
\bS \coloneqq& \left\{ \vect{p} = (p_0:p_1:p_2) \in \PP^2_A : \tilde P(\vect{p}) = 0\right\} \\
\isom& \left\{ \left( (b_{11},\dots,b_{33}), (p_0:p_1:p_2) \right) \in \Aff_\R^{6} \times \PP_\R^2 : \tilde P(\vect{p}) = 0\right\}.
\end{align*}
We call~$\bS$ the \defi{slowness bundle}.
Let $\iota\colon \bS \into\Aff_\R^{6} \times \PP_\R^2$ be the inclusion map. Composing~$\iota$ with the projection $\pi_1\colon\Aff_\R^{6} \times \PP_\R^2 \to \Aff_\R^{6}$ gives a fibration
\[
f \coloneqq \pi_1\circ\iota \colon \bS \to \Aff_\R^{6}
\]
that we call the \defi{slowness curve fibration}.  For a point $\vect{b} = (b_{11},\dots,b_{33}) \in \Aff^6(\R) = \R^6$, the fiber $f^{-1}(b)$ is the curve of degree~$4$ in the projective plane $\PP^2_\R$ obtained by specializing the parameters $b_{ij}$ in $\tilde P(\vect{p})$ according to the coordinates of~$\vect{b}$.

A theorem of Grothendieck known as ``generic geometric integrality''~\cite[Th\'eor\`eme 12.2.4(viii)]{EGAIV.3} allows us to conclude that the set of points $\vect{b} \in \R^6$ such that the fiber $f^{-1}(\vect{b})$ is irreducible, even over~$\C$, is an open subset for the Zariski topology of $\Aff_\R^6$.  This leaves two tasks for us: to show that the map~$f$ satisfies the hypotheses of~\cite[Th\'eor\`eme 12.2.4(viii)]{EGAIV.3} (i.e., that it is proper, flat, and of finite presentation), and that the open set in $\Aff_\R^6$ furnished by generic geometric irreducibility is not empty! For the latter, it suffices to produce a single choice of parameters~$b_{ij}$ such that the corresponding slowness curve is irreducible over~$\C$.  Since the target $\Aff_\R^6$ is an irreducible variety, the resulting open set is dense in the Zariski topology.

For the non-emptyness step, we use a standard number-theoretic strategy: reduction modulo a well-chosen prime.  To wit, we choose a slowness polynomial $\tilde P(\vect{p})$ with all $b_{ij} \in \Z$; to check it is irreducible over~$\C$, it suffices to show it is irreducible over a fixed algebraic closure~$\overline{\Q}$ of~$\Q$ (see~\cite[\href{https://stacks.math.columbia.edu/tag/020J}{Tag 020J}]{stacks-project}). Furthermore, a putative factorization would have to occur already over a finite Galois field extension $K\subset \overline{\Q}$ of~$\Q$, because all the coefficients involved in such a factorization would be algebraic numbers, and therefore have finite degree over~$\Q$.  Reducing the polynomial modulo a nonzero prime ideal $\frakp$ in the ring of integers $\calO_K$ of~$K$, by applying the unique ring homomorphism $\Z\to \calO_K/\frakp\calO_K =: \F_\frakp$ to its coefficients, we would see a factorization of $\tilde P(\vect{p})$ in the residue polynomial ring $\F_\frakp[\vect{p}]$, namely, the reduction of the factorization that occurs over~$K$. The finite field $\F_\frakp$ is an extension of the finite field~$\F_p$ with~$p$ elements, where $\frakp\cap \Z = (p)$. In Lemma~\ref{lem:IrreducibilityCriterion}, we show that if $\tilde P(\vect{p})$ is irreducible in the finite field~$\F_{p^d}$ of cardinality~$p^d$, where $d = \deg\big(\tilde P(\vect{p})\big)$, then it is also irreducible over the finite field~$\F_\frakp$, and hence is irreducible over~$K$, hence over~$\overline{\Q}$, hence over~$\C$. (This effectively shows that $K = \Q$.) What makes this strategy compelling is that~$\F_{p^d}$ is \emph{finite}, so checking whether $\tilde P(\vect{p})$ is irreducible in $\F_{p^d}[\vect{p}]$ is a finite, fast computation in any modern computer algebra system.

\begin{remark}
Readers versed in algebraic geometry might wonder if it might not be easier to use ``generic smoothness''~\cite[Corollary~III.10.7]{Hartshorne} to prove that a generic slowness polynomial is irreducible.  Unfortunately, in dimensions $n \notin \{2,4,8\}$, a slowness surface is \emph{always} singular \cite{I:degenerate}, and since $n = 3$ is the most interesting case from a physical point of view, we must avoid using ``generic smoothness.''
\end{remark}

%%%%%%%%%%%
\subsection{Irreducibility over \texorpdfstring{$\C$}{C} vs.\ connectedness over \texorpdfstring{$\R$}{R}}

It is possible for the set $S(\R) \subseteq \R^{2}$ of real points of a slowness curve~$S$ to be disconnected in the Euclidean topology, even if the algebraic variety~$S$, considered over the complex numbers, is connected in the Zariski topology.  For example, taking
\[
b_{11} = 10,\quad b_{12} = 2, \quad b_{13} = 3, \quad b_{22} = 12,\quad b_{23} = 5,\quad b_{33} = 20,
\]
we obtain the slowness curve (using coordinates $x = p_1$ and $y = p_2$):
\[
S : \quad 
1 - 30 x^2 + 191 x^4 - 16 x y + 88 x^3 y - 32 y^2 + 
 66 x^2 y^2 + 52 x y^3 + 215 y^4
= 0.
\]
This curve has two real connected components (see Figure~\ref{fig:2d-slowness-example}).
However, as a complex algebraic variety,~$S$ is irreducible, hence connected.  Its natural compactification in the complex projective plane is a smooth genus~$3$ complex curve, which is a $3$-holed connected $2$-dimensional real manifold.  

%%%%%%%%%%%
\subsection{Unique reconstruction}
\label{sec:uniquereconstruction}
Our second goal, unique reconstruction of generic stiffness tensors, has both a theoretical facet and a computational facet, which are in some sense independent.  Comparing the coefficients of~\eqref{eq:slownesspoly2D} and~\ref{eq:gen2Dpoly}, after dehomogenizing by setting $p_0 = 1$, we want to ideally solve the system of simultaneous equations
\begin{equation}
\begin{aligned}
\label{eq:2Dreconstructiondehom}
c_1 &= (b_{11}b_{33}-b_{13}^2), & c_5 &= 2(-b_{12}b_{23} + b_{13}b_{22}),\\
c_2 &= 2(b_{11}b_{23}-b_{12}b_{13}), & c_6 &= - 2(b_{13}+b_{23}),\\
c_3 &= (b_{11}b_{22}-b_{12}^2-2b_{12}b_{33}+2b_{13}b_{23}), & c_7 &= (b_{22}b_{33} - b_{23}^2),\\
c_4 &= - (b_{11}+b_{33}), & c_8 &= - (b_{22}+b_{33}).
\end{aligned}
\end{equation}
That is, given constants $c_1,\dots,c_8$, we would like to determine all $6$-tuples $(b_{11},\dots,b_{33})$ that satisfy~\eqref{eq:2Dreconstructiondehom}.
To this end, we homogenize the system in a slightly different way than before with a new variable $r$, so that all the right hand sides are homogeneous polynomials of degree $2$:
\begin{equation}
\begin{aligned}
\label{eq:2Dreconstruction}
\tilde c_1 &= (b_{11}b_{33}-b_{13}^2), & \tilde c_5 &= 2(-b_{12}b_{23} + b_{13}b_{22}),\\
\tilde c_2 &= 2(b_{11}b_{23}-b_{12}b_{13}), & \tilde c_6 &= - 2r(b_{13}+b_{23}),\\
\tilde c_3 &= (b_{11}b_{22}-b_{12}^2-2b_{12}b_{33}+2b_{13}b_{23}), & \tilde c_7 &= (b_{22}b_{33} - b_{23}^2),\\
\tilde c_4 &= - r(b_{11}+b_{33}), & \tilde c_8 &= - r(b_{22}+b_{33}).%,\\
\end{aligned}
\end{equation}
This homogenization allows us to define a \emph{rational map} of complex projective spaces
\begin{equation}
\label{eq:ratmap2D}
\begin{split}
g\colon \PP^6_{(b_{11}:\,\cdots\,:b_{33}:r)} &\dasharrow \PP^7_{(\tilde c_1:\,\cdots\,:\tilde c_8)}, \\
(b_{11}:\,\cdots\,:b_{33}:r) &\mapsto \left(b_{11}b_{33} - b_{13}^2:\,\cdots\,:-r(b_{22} + b_{33})\right).
\end{split}
\end{equation}
We want to show that a general nonempty fiber of this map consists of exactly one point; this would mean that among all tuples $(\tilde c_1,\dots,\tilde c_8)$ that are possibly the coefficients of a slowness polynomial, most tuples arise from exactly one stiffness tensor.

The map $g$ is strongly related to the map $\spolymap\colon \stiffnesstensor{2}\to\R[p_1,p_2]$ from \S\ref{ss:algebraicview}, but there are differences:
\smallskip
\begin{enumerate}[leftmargin=*]
    \item The indexing of the stiffness tensor elements is different (the map $g$ uses Voigt notation).
    \smallskip
    \item The map $\spolymap$ takes on polynomial values, namely, $F(\vect{a}) = P_\vect{a}(\vect{p})$, whereas the map $g$ outputs the coefficients of $\tilde P(\vect{p})$, which is a homogenized version of $P_\vect{a}(\vect{p})$.
    \smallskip
    \item The map $g$ is homogenized to map between projective spaces.
\end{enumerate}
\smallskip
Since $g$ contains the same information as $\spolymap$ modulo these identifications, to prove Theorem~\ref{thm:polynomial-to-tensor}(\ref{item:2D}), it suffices to prove that the map $g$ is generically injective. \\

The map $g$ is not defined at points where the right hand sides of~\eqref{eq:2Dreconstruction} simultaneously vanish. This is what is represented by the dashed arrow. Call this locus $\Pi \subset \PP^6$. The algebro-geometric operation of \emph{blowing up} $\PP^6$ at $\Pi$ gives a scheme $X \coloneqq \Bl_\Pi(\PP^6)$ together with a projection map $X \to \PP^6$ that resolves the indeterminacy locus of $g$, in the sense that the composition $X \to \PP^6 \dasharrow \PP^7$ can be extended to a full morphism $h \colon X \to \PP^7$ such that the triangle
\[
\xymatrix{
X \ar[d] \ar[dr]^{h} & \\
\PP^6 \ar@{-->}[r]^g & \PP^7
}
\]
commutes.

In the first proof of Theorem~\ref{thm:polynomial-to-tensor}(\ref{item:2D}), we show that the rational map $g\colon \PP^6 \dasharrow \PP^7$ is \emph{birational onto its image}. This means there exists a Zariski-dense open subset $U\subset \PP^6$, where $g$ is defined, such that $U \isom g(U)$ as algebraic varieties.  This is enough to show generic unique reconstruction of stiffness tensors.

Let~$Y$ be the image of~$h$.  In the second proof of Theorem~\ref{thm:polynomial-to-tensor}(\ref{item:2D}), we show there is a Zariski open subset of~$Y$ whose fibers are zero-dimensional, using upper semi-continuity of the fiber dimension function for the surjective map $h \colon X \to Y$.  Then, using upper semi-continuity of the degree function for finite morphisms, we show there is a Zariski open subset of~$Y$ whose fibers consist of precisely one point. 

As a bonus, we can use an effective version of Chevalley's Theorem~\cite{Barakat}, implemented in the package {\tt ZariskiFrames}~\cite{ZariskiFrames} to compute equations for the (set-theoretic) image of~$g$.  This, for example, explains how we arrived at the constraint~\eqref{eq:oneconditiononcs} for the tuple $(c_1,\dots,c_8)$. Tuples of coefficients $(c_1,\dots,c_8)$ in the image of~$g$ are said give rise to \defi{admissible slowness polynomials}. 

%%%%%%%%%%%
\subsection{In practice: Use Gr\"obner bases}
\label{ss:2DGB}
As a computational matter, given a specific tuple $(c_1,\dots,c_8)$ stemming from an admissible slowness polynomial, reconstructing its   stiffness tensor can be done essentially instantaneously using Gr\"obner bases. 
We work over the field $\Q$ so that we can use any one of several computational algebra systems with Gr\"obner bases packages, e.g., {\tt magma}, Macaulay2, Singular, Maple, or Sage. Thanks to Buchberger's criterion (see~\cite[\S2.6]{CLO}), it is possible to check the result of our calculations by hand, albeit laboriously.

A (reduced) Gr\"obner basis for an ideal $I \subset \Q[b_{11},b_{12},b_{13},b_{22},b_{23},b_{33}]$ under the lexicographic ordering $b_{11} > \cdots > b_{33}$ is a basis for~$I$ whose leading terms  generate the ideal consisting of the leading terms of \emph{all} polynomials in~$I$. In the event that there is exactly one tuple $(b_{11},b_{12},b_{13},b_{22},b_{23},b_{33}) \in \Q^6$ 
that satisfies the relations defining~$I$, this Gr\"obner basis will consist of precisely this tuple. For example, given the admissible slowness polynomial
\begin{equation*}
\tilde P(\vect{p}) = 
-3625p_1^4 + 1590p_1^3p_2 + 7129p_1^2p_2^2 - 50p_1^2p_0^2 + 8866p_1p_2^3 
+ 304p_1p_2p_0^2 - 8049p_2^4 - 14p_2^2p_0^2 + p_0^4
\end{equation*}
the following short piece of {\tt magma} code~\cite{magma} reconstructs the components of the stiffness tensor:
\begin{verbatim}
P<b11,b12,b13,b22,b23,b33> := PolynomialRing(Rationals(),6);
relations := [ 
         -3625 - (b11*b33 - b13^2), 
         1590 - 2*(b11*b23 - b13*b12),
         7129 - (b11*b22 + 2*b13*b23 - b12^2 - 2*b33*b12),
         -50 + (b11 + b33),
         8866 - 2*(b13*b22 - b23*b12),
         304 + 2*(b13 + b23),
         -8049 - (b33*b22 - b23^2),
         -14 + (b33 + b22) ];
I := ideal<P | relations>;
GroebnerBasis(I);
\end{verbatim}
The computation takes less than a millisecond, and returns
\begin{verbatim}
[
    b11 - 20,
    b12 - 39,
    b13 + 65,
    b22 + 16,
    b23 + 87,
    b33 - 30
]
\end{verbatim}
indicating the parameters $(b_{11},b_{12},b_{13},b_{22},b_{23},b_{33}) = (20,39,-65,-16,-87,30)$ of the \emph{only} stiffness tensor that gives rise to the specific polynomial $\tilde P(\vect{p})$ above, as the reader can check.

As mentioned above, this kind of Gr\"obner basis computation is independent of the theoretical result asserting that a generic slowness polynomial arises from a unique stiffness tensor.
In fact, these results complement each other nicely: Theorem~\ref{thm:polynomial-to-tensor}(\ref{item:2D}) implies that a Gr\"obner basis computation will succeed when applied to a generic slowness polynomial when $n=2$. Conjecture~\ref{conj:uniquereconstruction} would imply the same when $n=3$.

%%%%%%%%%%%%%%%%%%%%%%%%%%%%%%%%
\section{The two-layer model}
\label{sec:2-layer-def}

%%%%%%%%%%%
\subsection{Nested domains and stiffness tensors}

We say that $\omega$, $\Omega\subset\R^n$ are \defi{nested domains} if they are smooth, strictly convex, and bounded domains such that $\bar\omega\subset\Omega$. Let $\vect{a}$ and $\vect{A}\in\posstiffnesstensor{n}$ be two stiffness tensors associated  to the regions $\omega$ and $\Omega\setminus\omega$, respectively. We call these tensors \defi{admissible nested stiffness tensors} if the following conditions hold:
\medskip
\begin{itemize}
\item [(A1)]
For both tensors the largest eigenvalue of the Christoffel matrix $\Gamma(\vect{p})$ is simple for all $\vect{p}\neq0$. We refer to the corresponding subset of the slowness surface as the \defi{qP-branch}.
\medskip
\item [(A2)]
The qP-branch of the slowness surface of~$\vect{a}$ is inside that of~$\vect{A}$. (In other words, the slowness surfaces $s$, $S\subset\R^n$ of~$\vect{a}$ and~$\vect{A}$ are the boundaries of nested domains in the sense defined above.)
\end{itemize}

The two domains and the stiffness tensors are illustrated in Figure~\ref{fig:inverse}.

If the stiffness tensors~$\vect{a}$ and~$\vect{A}$ are isotropic, then the nestedness condition above simply means that the qP wave speed of~$\vect{a}$ is strictly higher than that of~$\vect{A}$. If~$\omega$ and~$\Omega$ are concentric balls, then the condition is equivalent with the Herglotz condition interpreted in a distributional sense; cf.~\cite{dHIK:layered-rigidity}. The Herglotz condition is widely used in the theory of geometric inverse problems as a generalization of the condition that a Riemannian manifold with a boundary has no trapped geodesics.

The piecewise constant stiffness tensor field corresponding to a pair of nested domains $\omega,\Omega$ and admissible nested stiffness tensors $\vect{a},\vect{A}$ is the function $\Omega\to\posstiffnesstensor{n}$ taking the value~$\vect{a}$ in~$\omega$ and~$\vect{A}$ in~$\Omega\setminus\omega$.

%%%%%%%%%%%
\subsection{Admissible rays}
Intuitively speaking, among all physically realizable piecewise linear ray paths, admissible rays are geometrically convenient ray paths. We make no claims about the amplitudes of the corresponding waves, but we expect most admissible rays to have non-zero amplitudes. We will describe separately the behaviour where the stiffness tensor is smooth and the behaviour at the interfaces~$\partial\omega$ and~$\partial\Omega$. Admissible ray paths will be piecewise linear paths satisfying certain conditions.

Suppose first that the stiffness tensor $\vect{a}(x)\in\posstiffnesstensor{n}$ is smooth. For every $(x,\vect{p})\in T^*\Omega$ the Christoffel matrix $\Gamma_{\vect{a}}(x,\vect{p})$ has~$n$ positive eigenvalues, possibly with repetitions. For any $m\in\{1,\dots,n\}$, let $G_{\vect{a}}^m\subset T^*M$ denote the subset where the $m$-th eigenvalue $\lambda_{\vect{a}}^{m}(x,\vect{p})$ of $\Gamma_{\vect{a}}(x,\vect{p})$ is simple. In this set the eigenvalue defines a smooth Hamiltonian $H^m_{\vect{a}}(x,\vect{p})=\frac12[\lambda_{\vect{a}}^{m}(x,\vect{p})^2-1]$. An \defi{admissible ray path} is the projection of an integral curve of the Hamiltonian flow from~$T^*\Omega$ to the base~$\Omega$. (The cotangent vector on the fiber we refer to as the momentum.)

In our setting~$\vect{a}$ is constant, so these integral curves are straight lines with constant speed parametrization. The speed depends on direction and polarization (or the eigenvalue index~$m$ or the branch of the slowness surface --- these are all equivalent).

At an interface two ray paths meet. We set two conditions for the incoming and outgoing {\it paths}:
\begin{itemize}
\item[(P1)]
Neither path is tangent to the interface. (This is convenient but ultimately unimportant.)
\medskip
\item[(P2)]
The component of the momentum tangent to the interface is the same for both incoming and outgoing rays.
\end{itemize}
The two meeting rays can be on the same or opposite sides of the interface, corresponding to reflected and refracted rays, respectively. The polarization is free to change.

The outer boundary~$\partial\Omega$ is also an interface. There the rays may either terminate (``refract to/from outside~$\Omega$'') or be reflected back in.

An admissible ray is a piecewise linear path, and we refer to the linear segments as legs.

\begin{remark}
Our definition of an admissible ray path excludes degenerate polarizations (which correspond to singular points on the slowness surface) and rays travelling along an interface. In the proof of Theorem~\ref{thm:2-layers} it is irrelevant whether these are included; their exclusion is not used nor would their inclusion be an issue. Rays tangent to an interface are irrelevant in the same way, as are the rare cases where the reflection or transmission coefficient is zero despite there being a kinematically possible ray.
\end{remark}

%%%%%%%%%%%
\subsection{Data}

We consider two kinds of data: pure travel time data (to be denoted by~$\nakeddata$) and travel time data decorated with direction information (to be denoted by~$\decorateddata$).

The full data set corresponding to the four parameters $(\omega,\Omega,\vect{a},\vect{A})$ is the set
\begin{equation*}
\begin{split}
\decorateddata(\omega,\Omega,\vect{a},\vect{A})
&=
\{
(t,x,p,y,q);
\,
x,y\in \partial\Omega,
\,
\text{there is an admissible ray path}
\\&\qquad
\text{from $x$ to $y$ with initial momentum $p$,}
\\&\qquad
\text{final momentum $q$, and total length $t$}
\}
.
\end{split}
\end{equation*}
The pure travel time data set without directional information is
\begin{equation*}
\nakeddata(\omega,\Omega,\vect{a},\vect{A})
=
\{
(t,x,y);
(t,x,p,y,q)\in \decorateddata(\omega,\Omega,\vect{a},\vect{A})
\text{ for some $p\in T_x^*\bar\Omega$ and $q\in T_y^*\bar\Omega$}
\}
.
\end{equation*}
These two sets may be seen as subsets:
$\decorateddata(\omega,\Omega,\vect{a},\vect{A})\subset\R\times(\partial T^*\bar\Omega)^2$
and
$\nakeddata(\omega,\Omega,\vect{a},\vect{A})\subset\R\times(\partial\Omega)^2$.

%%%%%%%%%%%%%%%%%%%%%%%%%%%%%%%%
\section{Inverse problems proofs}
\label{sec:IP-proofs}

This section is devoted to the proof of the inverse problems results, Theorems~\ref{thm:generic-stiffness-uniqueness} and~\ref{thm:2-layers}. We will make use of Theorems~\ref{thm:GenericIrreducibility} and~\ref{thm:polynomial-to-tensor} (and Conjecture~\ref{conj:uniquereconstruction} for $n=3$); besides them, we only need very basic algebraic geometry.

%%%%%%%%%%%
\subsection{Proof of Theorem \ref{thm:generic-stiffness-uniqueness}}

The first result follows easily from our algebraic results, Theorems~\ref{thm:GenericIrreducibility} and~\ref{thm:polynomial-to-tensor} that we prove in~\S\ref{sec:AG-proofs}.

\begin{proof}[Proof of Theorem \ref{thm:generic-stiffness-uniqueness}]
Theorem~\ref{thm:GenericIrreducibility} implies that there is an open and dense (in the Zariski sense) set $W_1\subset\stiffnesstensor{n}$ so that the slowness polynomial~$P_\vect{a}$ is irreducible for all $a\in W_1$. The Zariski-closure of the relatively open (in the Euclidean sense) subset of the slowness surface is a subvariety of the slowness surface, and it is of full dimension. Due to irreducibility this closure is the whole slowness surface. Thus for $\vect{a}\in W_1$ a small subset of the slowness surface determines the whole slowness surface. This proves part~(\ref{item:thmA-surface}), so we move to part~(\ref{item:thmA-tensor}).

Theorem~\ref{thm:polynomial-to-tensor} ($n=2$) or Conjecture~\ref{conj:uniquereconstruction} ($n=3$) implies that there is an open and dense set $W_2\subset\stiffnesstensor{n}$ so that if $P_{\vect{a}_1}=P_{\vect{a}_2}$ for $\vect{a}_1,\vect{a}_2\in W_2$, then $\vect{a}_1=\vect{a}_2$.

The set $W\coloneqq W_1\cap W_2$ is open and dense (in the Zariski sense) in $\stiffnesstensor{n}$. If $\vect{a}_1,\vect{a}_2\in W$, then the corresponding slowness surfaces agreeing in a any open subset implies by part~(\ref{item:thmA-surface}) that the slowness surfaces agree. By the argument above this implies that the stiffness tensors agree as well: $\vect{a}_1=\vect{a}_2$.

The positivity property of the stiffness tensor was irrelevant for both parts. The claims remain true in the physically relevant open subset $\posstiffnesstensor{n}\subset\stiffnesstensor{n}$ by taking $U=W\cap\posstiffnesstensor{n}$.
\end{proof}

%%%%%%%%%%%
\subsection{Proof of Theorem \ref{thm:2-layers}}

This proof will rely on Theorem~\ref{thm:generic-stiffness-uniqueness} without the use any algebraic geometry. We will split the proof in three parts, proven separately below:

\begin{enumerate}
\item\label{step:2-layer-part-1}
If $\nakeddata(\omega_1,\Omega,\vect{a}_1,\vect{A}_1)
\approx
\nakeddata(\omega_2,\Omega,\vect{a}_2,\vect{A}_2)$, then $\vect{A}_1=\vect{A}_2$.
\medskip
\item\label{step:2-layer-part-2}
If $\nakeddata(\omega_1,\Omega,\vect{a}_1,\vect{A}_1)
\approx
\nakeddata(\omega_2,\Omega,\vect{a}_2,\vect{A}_2)$, then $\omega_1=\omega_2$.
\medskip
\item\label{step:2-layer-part-3}
If $\decorateddata(\omega_1,\Omega,\vect{a}_1,\vect{A}_1)
\approx
\decorateddata(\omega_2,\Omega,\vect{a}_2,\vect{A}_2)$, then $\vect{a}_1=\vect{a}_2$.
\end{enumerate}

\begin{figure}
    \centering
    \begin{overpic}[width=0.4\textwidth]{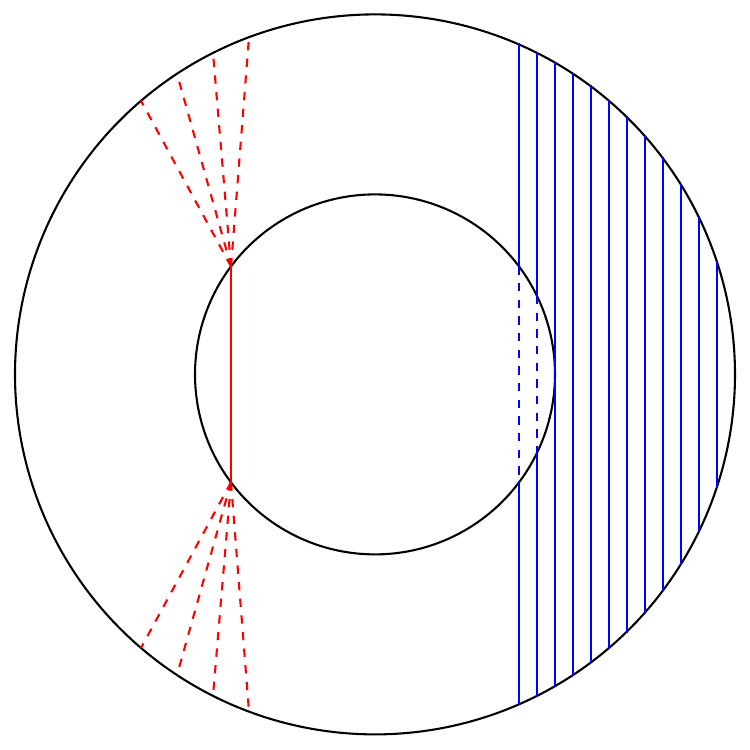}
    \put(47,69){$\partial\omega_i$}
    \put(47,100){$\partial\Omega$}
    \put(40,45){$\vect{a}$}
    \put(15,45){$\vect{A}$}
    \end{overpic}
    \caption{A cartoon of the proof of Theorem~\ref{thm:2-layers}.
    In the first step, we use short geodesics near the boundary, depicted as the \textcolor{blue}{blue line segments to the right}.
    In the second step, we vary this family of line segments until they hit~$\partial\omega_i$ and no longer exist, depicted as the \textcolor{blue}{dashed parts of the blue line segments}.
    In the third step we take two points on~$\partial\omega_i$ and find their qP distance, depicted as \textcolor{red}{the solid red line}, hitting them with all possible rays through the now known mantle, depicted as \textcolor{red}{the dashed red lines}.
    }
    \label{fig:inverse}
\end{figure}
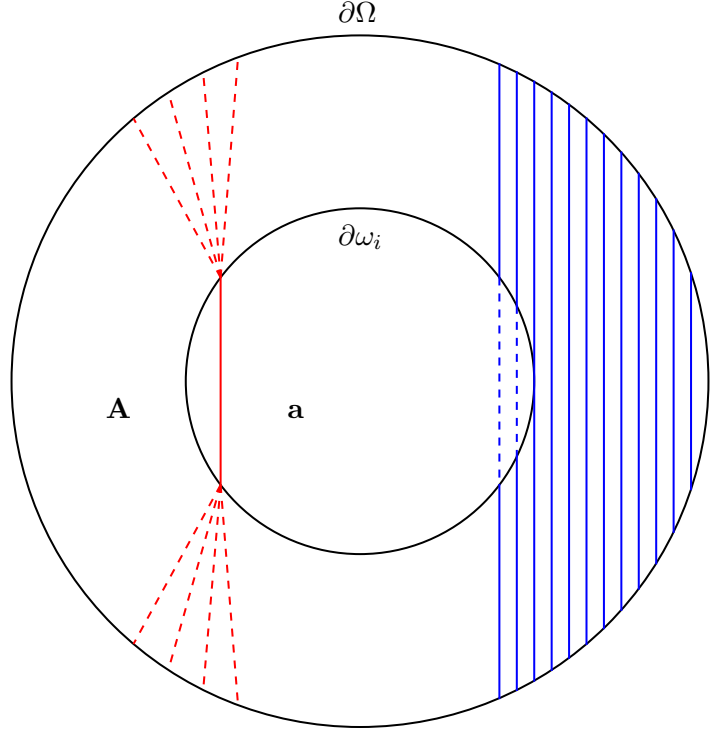

Roughly speaking, we will prove the first part by studying the travel times of nearby points, the second part by varying a line segment and detecting when it hits~$\partial\omega_i$, and the third part by peeling off the top layer to get a problem on~$\partial\omega$ that is similar to the first step. These parts are illustrated in Figure~\ref{fig:inverse}.

For any $x\in\partial\Omega$, let $\nu(x)$ denote the inward-pointing unit normal vector to the boundary~$\partial \Omega$.
Given any direction $d\in\Sphere^{n-1}$, we denote
\begin{equation*}
\partial_d\Omega
=
\{
x\in\partial\Omega
;
d\cdot\nu(x)>0
\}
.
\end{equation*}
This is the subset of the boundary where~$d$ points inwards. The boundary of this set, $\partial\partial_d\Omega\subset\partial\Omega$, is the set where~$d$ is tangent to the boundary.

Due to strict convexity of~$\Omega$ there is a unique $\tau(x,d)>0$ for every $x\in\partial_d\Omega$ so that $x+\tau(x,d)d\in\partial\Omega$. This is the distance through~$\Omega$ starting at~$x$ in the direction~$d$. For $x\in\partial\partial_d\Omega$ we set $\tau(x,d)=0$, and we do not define the function at all where~$d$ points outwards.

For both $i=1,2$, we denote
\begin{equation}
\label{eq:vi(x,d)}
v_i(x,d)
=
\frac
{
\tau(x,d)
}
{
\inf\{t>0;(t,x,x+\tau(x,d)d)\in\nakeddata(\omega_i,\Omega,\vect{a}_i,\vect{A}_i)\}
}
.
\end{equation}
This can be seen as a speed through~$\Omega$ starting from~$x$ and ending in the direction~$d$ from~$x$, but not knowing the initial and final directions of the minimizing ray or whether the ray has reflected from the interfaces $\partial\omega$ or $\partial\Omega$, or whether it has met $\partial\omega$ tangentially. This function~$v_i$ is a suitable form of data for the first two steps of the proof.

\begin{lemma}
\label{lma:fastest-route}
Let $n\geq2$. Let $\omega,\Omega\subset\R^n$ be nested domains. Take any $d\in\Sphere^{n-1}$ and $x\in\partial_d\Omega$. Denote $y\coloneqq x+\tau(x,d)d$.

Let $\vect{a},\vect{A}\in\posstiffnesstensor{n}$ be two stiffness tensors whose Christoffel matrices $\Gamma_{\vect{a}}(p)$ and $\Gamma_{\vect{A}}(p)$ have a simple largest eigenvalue for all $p\neq0$.

\begin{enumerate}
\item[a)]
If $\vect{a}=\vect{A}$ or $x+d\R\cap\bar\omega=\emptyset$, then the fastest admissible ray path between~$x$ and~$y$ is the qP polarized ray travelling along the straight line between the points.
\item[b)]
If~$\vect{a}$ and~$\vect{A}$ are admissible nested stiffness tensors and $x+d\R\cap\omega\neq\emptyset$, then the shortest travel time between~$x$ and~$y$ is strictly larger than it would be if~$\vect{a}$ were equal to~$\vect{A}$.
\end{enumerate}
\end{lemma}

\begin{proof}
a)
The qP slowness surface is strictly convex as observed in \cite{dHILS:BDF}, so the integral curves of the Hamiltonian flow do indeed minimize length. With a constant stiffness tensor this minimization property is global.

b)
The nestedness property of the qP branches of the slowness surfaces imply that all ray paths for the tensor~$\vect{a}$ are slower than those of~$\vect{A}$ in the same direction. Therefore for every admissible ray path that meets~$\omega$ the total travel time is strictly bigger than the length of that piecewise smooth curve measured in the qP Finsler geometry of~$\vect{A}$. Therefore the shortest travel time of an admissible ray path has to be longer than it would be if~$\vect{a}$ were changed to be equal to~$\vect{A}$.
\end{proof}

We will denote the data sets by $\decorateddata_i\coloneqq\decorateddata(\omega,\Omega,\vect{a}_i,\vect{A}_i)$ and~$\nakeddata_i$ similarly.

\begin{proof}[Proof of Theorem \ref{thm:2-layers}, part~\ref{step:2-layer-part-1}]
The set~$U$ is taken to be that provided by Theorem~\ref{thm:generic-stiffness-uniqueness}.

The functions $v_i(x,d)$ of~\eqref{eq:vi(x,d)} are defined in the subset
\begin{equation*}
\{
(x,d)
\in
\partial\Omega\times\Sphere^{n-1}
;
x\in\overline{\partial_d\Omega}
\}
\end{equation*}
and are continuous in a neighborhood~$\mathcal{U}$ of the subset $\tau^{-1}(0)$, corresponding to short geodesics that do not meet~$\bar\omega_i$.

The assumption
$\nakeddata_1
\approx
\nakeddata_2$
implies that the functions~$v_1$ and~$v_2$ agree in an open and dense subset of~$\mathcal{U}$, and by continuity they agree on all of~$\mathcal{U}$. Thus near the boundary we may work as if
$\nakeddata_1
=
\nakeddata_2$.

In fact, these functions~$v_i$ only depend on~$d$ in~$\mathcal{U}$. Fix any direction $d\in\Sphere^{n-1}$. By strict convexity of the nested domains~$\omega$ and~$\Omega$, there is a neighborhood $Y\subset\partial\Omega$ of $\partial\partial_d\Omega$ so that for all $x\in\hat Y\coloneqq Y\cap\partial_d\Omega$ the ray starting at~$x$ in the direction~$d$ does not meet $\bar\omega_1\cup\bar\omega_2$. (We remind the reader that the direction~$d$ is tangent to the boundary precisely in the set~$\partial\partial_d\Omega$. Therefore in a small neighborhood of this set the line segments in the direction~$d$ through~$\Omega$ are short.) By Lemma~\ref{lma:fastest-route} a qP polarized ray travelling from~$x$ in the direction~$d$ minimizes the travel time between~$x$ and $x+\tau(x,d)d$.

This implies that both functions $v_i(\dummy,d)$ are constant in~$\hat Y$. By the assumption of the agreement of the data~$\nakeddata$, these two functions agree. Let us denote to shared constant value by $v(d)$. Therefore the two models give rise to the same surfaces
\begin{equation*}
S^*
=
\{
v(d)d
;
d\in\Sphere^{n-1}
\}
.
\end{equation*}
The surface $S^*$ is the strictly convex unit sphere of the Finsler geometry corresponding to the qP polarized waves; cf.~\cite[\S~2]{dHILS:BDF}. By taking the Legendre transform, the set~$S^*$ determines the dual sphere $S\subset\R^n$ in the usual sense of dual norms. This cosphere~$S$ is exactly the qP branch of the slowness surface.

By assumption each~$\vect{A}_i$ is in~$U$, the open and dense set provided by Theorem~\ref{thm:generic-stiffness-uniqueness}. Therefore this branch of the slowness surface determines the stiffness tensor, and so $\vect{A}_1=\vect{A}_2$.
\end{proof}

\begin{proof}[Proof of Theorem \ref{thm:2-layers}, part~\ref{step:2-layer-part-2}]
We denote $\vect{A}\coloneqq \vect{A}_1=\vect{A}_2$.

Again, fix any direction $d\in\Sphere^{n-1}$. Let $Y_i^d\subset\partial_d\Omega$ be the subset where $v_i(x,d)$ takes the constant value $v(d)$; cf. part~\ref{step:2-layer-part-1} of the proof. Let us denote $\mathcal{U}_i'=\{(x,d);x\in Y_i^d\}$. As the data is defined as subsets of (the real axis and two copies of) the set $\partial\Omega\times\Sphere^{n-1}$, it follows from approximate equality of the data (the assumption $\nakeddata_1\approx\nakeddata_2$) that
\begin{equation*}
\overline{\mathcal{U}_1'}
=
\overline{\mathcal{U}_2'}
\eqqcolon
\overline{U'}
.
\end{equation*}
We will use this set to describe the inner domains~$\omega_i$.

It follows from Lemma~\ref{lma:fastest-route} and the definition of~$v_i$ as a directed travel time that if the line $x+d\R$ meets~$\omega_i$, then $x\notin Y_i^d$, and that if it does not meet~$\bar\omega_i$, then $x\in Y_i^d$. The line $x+d\R$ is tangent to~$\partial\omega_i$ if and only if $x\in\partial Y_i^d$. We thus know which lines meet~$\omega_i$, and can write the smaller domain as
\begin{equation*}
\omega_i
=
\Omega
\setminus
\overline{
\bigcup_{(x,d)\in\mathcal{U}'_i}
(x+d\R)
}
=
\Omega
\setminus
\bigcup_{(x,d)\in\overline{\mathcal{U}'}}
(x+d\R)
.
\end{equation*}
Therefore $\omega_1=\omega_2$ as claimed.
\end{proof}

We can rephrase the proof above loosely as follows. We may think that $Y_1^d=Y_2^d$ (although this was not assumed to hold perfectly and for all~$d$) and say that the two strictly convex and smooth domains~$\omega_1$ and~$\omega_2$ have the same tangent lines so they are equal.

\begin{proof}[Proof of Theorem \ref{thm:2-layers}, part~\ref{step:2-layer-part-3}]
As in the previous proofs, we can essentially replace the assumption $\decorateddata_1\approx\decorateddata_2$ with the stronger one $\decorateddata_1=\decorateddata_2$ because we are using open subsets of the data sets rather than relying on rare features. We omit the details in this instance for clarity.

By the previous parts of the theorem, we now know that $\vect{A}_1=\vect{A}_2\eqqcolon \vect{A}$ and $\omega_1=\omega_2\eqqcolon\omega$. It remains to show that $\vect{a}_1=\vect{a}_2$. As in part~\ref{step:2-layer-part-1}, it suffices to prove that some non-empty open subsets of the qP branches of the two slowness surfaces agree.

Each point of $\partial\Omega\times\Sphere^{n-1}$ defines a ray starting at the given point on~$\partial\Omega$ in the given direction in~$\Sphere^{n-1}$. For any $x\in\Omega$, there is a subset $F_x\subset\partial\Omega\times\Sphere^{n-1}$ so that the corresponding rays meet~$x$. The set~$F_x$ may be thought of as the graph of the unit vector field on~$\partial\Omega$ pointing towards~$x$. Let $F_x'\subset F_x$ be the subset corresponding to rays that do not meet~$\omega$ before~$x$.

Let $f_{\vect{A}}\colon\R^n\to[0,\infty)$ be the smooth and strictly convex norm whose unit sphere is the qP branch of the slowness surface corresponding to the stiffness tensor~$\vect{A}$. Let $f^*_{\vect{A}}$ be the dual norm and let $\phi_{\vect{A}}\colon\R^n\to\R^n$ be the norm-preserving and homogeneous (but possibly non-linear) Legendre transformation satisfying $f(p)^2=\ip{\phi_{\vect{A}}(p)}{p}$. For a direction $v\in\Sphere^{n-1}$, let us denote $Q_{\vect{A}}(v)=\phi_{\vect{A}}^{-1}(v/f_{\vect{A}}^*(v))$. In words, $Q_{\vect{A}}(v)$ is the momentum corresponding to the qP polarized wave travelling in the direction~$v$ in a material given by~$\vect{A}$. The Legendre transform is depicted in Figure~\ref{fig:2d-slowness-example}, where points (or arrows) on the cotangent space correspond to arrows (or points, respectively) on the tangent space.

Let us then define
\begin{equation*}
F_x''
=
\{
(z,Q_{\vect{A}}(v));
(z,v)\in F_x'
\}
.
\end{equation*}
This is the set of qP momenta (instead of directions) on the boundary so that the corresponding rays meet~$x$ without hitting~$\bar\omega$ first.

For each $(z,p)\in F_x''$ the travel time (according to the Hamiltonian flow) of the qP wave from~$z$ to~$x$ is $f_{\vect{A}}^*(x-z)$.

Now let $x,y\in\partial\omega$ be two distinct points and define
\begin{equation*}
\begin{split}
T_i(x,y)
&=
\inf\{
t
-f_{\vect{A}}^*(x-z)
-f_{\vect{A}}^*(y-w)
;
\\&\qquad
(t,z,p,w,-q)\in \decorateddata_i
\\&\qquad
\text{and } (z,p)\in F_x''
\text{ and } (w,q)\in F_y''
\}
.
\end{split}
\end{equation*}
Each ray considered here starts with a qP polarized leg from a point $z\in\partial\Omega$ to $x\in\partial\omega$ and ends in a similar leg from~$y$ to~$w$. As the travel times of the first and last legs are removed from the total travel time, our $T(x,y)$ is the shortest travel time between the points~$x$ and~$y$ with an admissible ray path.

Because the qP branches of the slowness surfaces of~$\vect{a}_i$ and~$\vect{A}$ are nested by assumption, all momenta are available for the segments of the ray path in~$\omega$ starting at~$x$ and~$y$.

All the travel times and the geometry between~$z$ and~$x$ and also between~$y$ and~$w$ are the same between the two models by the previous two steps, and the only remaining dependence on~$i$ is what happens between~$x$ and~$y$.

We claim that when~$x$ and~$y$ are sufficiently close to each other,
\begin{equation}
\label{eq:T=f*}
T_i(x,y)
=
f_{\vect{a}_i}^*(x-y).
\end{equation}
This means that the shortest admissible ray path between~$x$ and~$y$ is the direct qP ray within~$\omega$. This is seen as follows: If a ray path has a leg in the outer layer $\Omega\setminus\bar{\omega}$ between $x$ and~$y$ (which may well happen, as we do not a priori know the geometry of the rays we are looking at), then by strict convexity of~$\omega$ this leg must come all the way to the outer boundary~$\partial\Omega$. If~$x$ and~$y$ are so close to each other that $f_{\vect{a}_i}^*(x-y)$ is less than the $f_{\vect{A}}^*$-distance between $\partial\omega$ and~$\partial\Omega$, then any leg joining $\partial\omega$ to~$\partial\Omega$ takes a longer time than the straightest option through~$\omega$, despite the waves being slower in~$\omega$ than in~$\Omega\setminus\bar\omega$. Within~$\bar\omega$ the shortest travel time is clearly achieved by going in a straight line with the fastest polarization; cf. the Lemma~\ref{lma:fastest-route}.

If we fix $x\in\partial\omega$, we have found that equation~\eqref{eq:T=f*} holds for all~$y$ in a small punctured neighborhood of~$x$ on~$\partial\omega$ for both $i=1,2$. Because $\decorateddata_1=\decorateddata_2$ implies $T_1(x,y)=T_2(x,y)$, we have found that there is an open set $Y_x\subset\partial\omega$ so that
\begin{equation*}
f_{\vect{a}_1}^*(x-y)
=
f_{\vect{a}_2}^*(x-y)
\end{equation*}
for all $y\in Y_x$.
By strict convexity of~$\partial\omega$ the set~$Y_x$ contains an open set of directions, so the unit spheres of $f_{\vect{a}_1}^*$ and $f_{\vect{a}_2}^*$ agree on an open set. The same thus holds for $f_{\vect{a}_1}$ and $f_{\vect{a}_2}$ as well, and the claim $\vect{a}_1=\vect{a}_2$ follows from Theorem~\ref{thm:generic-stiffness-uniqueness}.
\end{proof}

%%%%%%%%%%%%%%%%%%%%%%%%%%%%%%%%
\section{The algebraic geometry of families of slowness surfaces}
\label{sec:AG-proofs}

This section contains the technical algebro-geometric arguments needed to prove Theorems~\ref{thm:GenericIrreducibility} and~\ref{thm:polynomial-to-tensor}; it demands more expertise from the reader than \S\ref{sec:alg-principles}. Our arguments use only material typically covered in a first course on scheme-theoretic algebraic geometry. Standard references for this material include~\cites{Hartshorne,Liu,Vakil,GWAG}. We provide copious references to specific propositions and theorems to help orient readers less familiar with schemes.

%%%%%%%%%%%
\subsection{Independent components of a stiffness tensor: Voigt notation}
\label{sec:tensor-alg}

The major and minor symmetries of a (reduced) stiffness tensor allow for a simplification of notation that eliminates clutter, following \defi{Voigt}. In dimension~$2$, one replaces pairs of indices~$ij$ by a single index~$k$ according to the rule
\begin{equation}
\label{eq:voigt2D}
11\leadsto1,\quad 22\leadsto2,\quad 12\leadsto3.
\end{equation}
To avoid confusion, when we contract indices following this convention, we also replace the letter~$a$ with the letter~$b$: for example, the reduced stiffness tensor component $a_{1112} = a_{(11)(12)}$ is replaced by~$b_{13}$. 

In dimension~$3$ one replaces pairs of indices~$ij$ by a single index~$k$ according to the rule
\begin{equation}
\label{eq:voigt3D}
11\leadsto1,\quad
22\leadsto2,\quad
33\leadsto3,\quad
23\leadsto4,\quad
13\leadsto5,\quad
12\leadsto6.
\end{equation}
Thus, for example, the reduced stiffness tensor component $a_{2312} = a_{(23)(12)}$ is replaced by~$b_{46}$. 

Next, we count the number of independent parameters of the form~$a_{ijkl}$, or equivalently, the number of independent parameters of the form~$b_{rs}$, once we take the symmetries~\eqref{eq:symmetry} into account.  The set of distinct~$a_{ijkl}$ is in bijection with a set of unordered pairs of unordered pairs of indices $\{1,\dots,n\}$: more precisely, a set whose elements have the form $a_{(ij)(kl)}$, where the indices belong to $\{1,\dots,n\}$ and one can freely commute the indices within a pair of parentheses or commute the pairs, but one cannot freely move indices from one pair to another.
The number of unordered pairs of indices $1,\dots,n$ is $\pairs(n) \coloneqq\frac12n(n+1)$.
Therefore the number of independent components of a stiffness tensor is
\[
\pairs(\pairs(n))=\frac18n(n+1)(n^2+n+2).
\]
When $n=2$, we obtain $\pairs(\pairs(2)) = 6$, which matches our work in \S\ref{sec:alg-principles}, where we saw the six independent parameters $b_{11}$, $b_{12}$, $b_{13}$, $b_{22}$, $b_{23}$, and $b_{33}$. In dimension $n = 3$, there are $\pairs(\pairs(3)) = 21$ independent stiffness tensor components.

%%%%%%%%%%%
\subsection{Algebro-geometric set-up}
\label{sec:alg-setup}

In this section, we omit the positivity condition that a stiffness tensor satisfies (Definition~\ref{def:stiffness-tensor}), in order to import ideas from the scheme-theoretic formulation of algebraic geometry, following Grothendieck.  

\subsubsection{The slowness polynomial}

Let~$A$ be a finitely generated $\Q$-algebra, and let $R \coloneqq A[p_0,\dots,p_{n}]$ be a polynomial ring in $n+1$ variables with coefficients in~$A$. We view the Christoffel matrix~\eqref{eq:ChristoffelEntry} as a symmetric $n\times n$ matrix $\Gamma(\vect{p})$ with entries in $R$, whose $il$-th entry is
\[
\Gamma(\vect{p})_{il}
=
\sum_{1\leq j,k \leq n}
a_{ijkl}p_jp_k,
\]
and where the parameters  $a_{ijkl} \in A$ are subject to the symmetry relations~\eqref{eq:symmetry}. Denoting by~$I_n$ the $n\times n$ identity matrix over~$A$, the \defi{slowness polynomial} $\spoly{\vect{a}}(\vect{p}) \in R$ is
\begin{equation*}
\spoly{\vect{a}}(\vect{p}) \coloneqq \det(\Gamma(\vect{p})-I_n).
\end{equation*}
This is a polynomial of total degree $d=2n$ in $p_1,\dots,p_{n}$.  The \defi{homogenized slowness polynomial} $\tilde P(\vect{p}) \in R$ is obtained by setting
\begin{equation*}
\tilde P(\vect{p})
\coloneqq
\det(\Gamma(\vect{p}) - p_0^2I_n).
\end{equation*}
The \defi{completed slowness hypersurface}~$\bS$ is the algebraic hypersurface in the projective space~$\PP_A^{n}$ where~$\tilde P$ vanishes.  More precisely, the quotient ring homomorphism 
\[
A[p_0,\dots,p_{n}] \to A[p_0,\dots,p_{n}]/(\tilde P)
\]
describes a closed embedding $\bS \into \PP_A^n$ via the~$\Proj$ construction (see, for example,~\cite[\S II.2 and Exercise II.3.12]{Hartshorne}).

From now on, we specialize to the case 
\[
A = \Q[a_{ijkl} : 1\leq i,j,k,l \leq n],
\]
where the~$a_{ijkl}$ are indeterminates subject to the symmetry relations~\eqref{eq:symmetry}. By~\S\ref{sec:tensor-alg}, the ring~$A$ is a free $\Q$-algebra on $
m \coloneqq \psi(\psi(n))= \frac18n(n+1)(n^2+n+2)$ generators. When we consider orthorhombic and monoclinic stiffness tensors, we will further specialize $A$ by removing appropriate variables, using the conventions in Notation~\ref{nt:orthomono}.

\begin{example} 
\label{ex:n=2}
Let $n = 2$.  Then $A = \Q[a_{ijkl} : 1\leq i,j,k,l \leq 2]$, and there are only $\psi(\psi(2)) = 6$ distinct $a_{ijkl}$'s, which we relabel $b_{11}$, $b_{12}$, $b_{13}$, $b_{22}$, $b_{23}$, and $b_{33}$ using Voigt notation~\eqref{eq:voigt2D} as we did in \S\ref{sec:alg-principles}. Thus,~$A$ is the polynomial ring $\Q[b_{11},b_{12},b_{13},b_{22},b_{23},b_{33}]$, and the Christoffel matrix $\Gamma(\vect{p})$ is given as in~\eqref{eq:2DChristoffelmatrix}, with associated homogenized slowness polynomial $\tilde P(\vect{p})$ as in~\eqref{eq:slownesspoly2D}. The associated completed slowness hypersurface~$\bS$ is a quartic curve on~$\PP^2_A$ defined by the condition $\tilde P(\vect{p}) = 0$.
\end{example}

\begin{example} 
\label{ex:n=3}
Let $n = 3$.  Then $A = \Q[a_{ijkl} : 1\leq i,j,k,l \leq 3]$, and there are only $\psi(\psi(3)) = 21$ distinct $a_{ijkl}$'s, which we relabel $\{b_{ij} : 1 \leq i \leq j \leq 6\}$ using Voigt notation~\eqref{eq:voigt3D}. Thus,~$A$ is the polynomial ring $\Q[b_{11},\dots,b_{66}]$ in $21$ variables, and the associated homogenized slowness polynomial $\tilde P(\vect{p})$ has $50$ terms. The associated completed slowness hypersurface~$\bS$ is a sextic surface on~$\PP^3_A$ defined by the condition $\tilde P(\vect{p}) = 0$.
\end{example}

\begin{example}
    \label{ex:n=3ortho}
    Let $n = 3$. In the context of orthorhombic stiffness tensors there are only $9$ distinct parameters among the 81 elements of $\{a_{ijkl}, 1 \leq i, j, k, l \leq 3\}$ (see Notation~\ref{nt:orthomono}). Using Voigt's notation, these are
    \[
    b_{11},\ b_{12},\ b_{13},\ b_{22},\ b_{23},\ b_{33},\ b_{44},\ b_{55},\quad\text{and}\quad b_{66}.
    \]
    Thus, we set $A = \Q[b_{11}, b_{12}, b_{13}, b_{22}, b_{23}, b_{33}, b_{44}, b_{55},b_{66}]$. The Christoffel matrix of an orthorhombic stiffness tensor is
\[
\Gamma(\vect{p}) =
\begin{pmatrix}
b_{11}p_1^2 + b_{66}p_2^2 + b_{55}p_3^2  & (b_{12} + b_{66})p_1p_2 &  (b_{13} + b_{55})p_1p_3 \\
(b_{12} + b_{66})p_1p_2  & b_{66}p_1^2 + b_{22}p_2^2 + b_{44}p_3^2  & (b_{23} + b_{44})p_2p_3 \\
(b_{13} + b_{55})p_1p_3  & (b_{23} + b_{44})p_2p_3 & b_{55}p_1^2 + b_{44}p_2^2 + b_{33}p_3^2
\end{pmatrix}
\]
The homogenized slowness polynomial of such a tensor has the form
\begin{equation}
\label{eq:genericorthorhombic}
\begin{split}
\tilde P(\vect{p}) &= \det(\Gamma(\vect{p} - p_0^2I_3) \\
&=
    c_{1}p_1^6 + c_{2}p_1^4p_2^2 + c_{3}p_1^4p_3^2 + c_{4}p_1^4p_0^2 + c_{5}p_1^2p_2^4 + c_{6}p_1^2p_2^2p_3^2 + c_{7}p_1^2p_2^2p_0^2 \\ &\quad+ c_{8}p_1^2p_3^4
    + c_{9}p_1^2p_3^2p_0^2 + c_{10}p_1^2p_0^4 + c_{11}p_2^6 + c_{12}p_2^4p_3^2 + c_{13}p_2^4p_0^2\\ &\quad+ c_{14}p_2^2p_3^4 + c_{15}p_2^2p_3^2p_0^2 + c_{16}p_2^2p_0^4 + c_{17}p_3^6 + c_{18}p_3^4p_0^2 + c_{19}p_3^2p_0^4 + c_{20}p_0^6,
\end{split}    
\end{equation}
where, for example, we have
\begin{equation}
\label{eq:coeffexample}
    c_{7} = -b_{11}b_{22} - b_{11}b_{44} + b_{12}^2 + 2b_{12}b_{66} - b_{22}b_{55} - b_{44}b_{66} - b_{55}b_{66}.
\end{equation}
The associated completed slowness hypersurface~$\bS$ is a sextic surface on~$\PP^3_A$ defined by the condition $\tilde P(\vect{p}) = 0$.
\end{example}

\subsubsection{The slowness fibration}
\label{sss:SlownessFibration}
Generalizing Examples~\ref{ex:n=2}, \ref{ex:n=3}~and \ref{ex:n=3ortho} the homogenized slowness polynomial can be viewed as a homogeneous polynomial of degree $2n$ in the graded ring $A[p_0,\dots,p_n]$, where $A$ is itself a polynomial ring in, say, $m$ variables. From this perspective, the completed slowness hypersurface may be viewed as a hypersurface in the product of an affine space and a projective space:
\[
\bS = \left\{ \tilde P(\vect{b}) = 0\right\} \subset \PP^n_A \isom \Aff_\Q^{m} \times \PP_\Q^n = \Spec A \times_{\Spec \Q} \Proj \Q[v_0,\dots,v_n].
\]
We call~$\bS$ the \defi{slowness bundle}, and denote this closed immersion by $\iota\colon \bS \into\Aff_\Q^{m} \times \PP_\Q^n$. Composing~$\iota$ with the projection $\pi_1\colon\Aff_\Q^{m} \times \PP_\Q^n \to \Aff_\Q^{m}$ onto the first factor gives a fibration
\[
f \coloneqq \pi_1\circ\iota \colon \bS \to \Aff_\Q^{m}
\]
that we call the \defi{slowness surface fibration}.  The fiber $f^{-1}(\vect{b})$ of~$f$ above a rational point $\vect{b} \in \Aff^{m}(\Q) = \Q^{m}$ is the hypersurface of degree $2n$ in~$\PP^n_\Q$ obtained by specializing the parameters~$b_{ij}$ according to the coordinates of~$\vect{b}$.

For a field extension $K/\Q$, we write $f_{K}\colon S_K \to \Aff_K^{m}$ for the slowness surface fibration obtained as above after replacing~$\Q$ by~$K$ everywhere.
This is known in algebraic geometry as the ``base-extension of the morphism~$f$ by the map $\Spec K \to \Spec \Q$''. We are mostly interested in the cases $K = \R$ and $K = \C$. We call $f_\C \colon \bS_\C \to \Aff_\C^{m}$ the \defi{complexified slowness surface fibration}.

%%%%%%%%%%%
\subsection{Key results}
\label{sec:results}

The precise result underpinning Theorem~\ref{thm:GenericIrreducibility} is the following.
\begin{theorem}
\label{thm:precisegeomirreduc}
The set 
\[
\Irr(f) \coloneqq \{\vect{b} \in \Aff^{m}_\R : f_{\C}^{-1}(\vect{b}) \textup{ is an irreducible hypersurface}\}
\]
is Zariski-open in $\Aff^{m}_\R$. Consequently, it is either empty, or it is the complement of a finite union of algebraic varieties, each of dimension $\leq m-1$.
\end{theorem}

Theorem~\ref{thm:precisegeomirreduc} follows from the following result, due to Grothendieck, which uses the full power of scheme theory.

\begin{theorem}[Generic Geometric Irreducibility]
\label{thm:GenericIrreducibilityEGA}
Let $f\colon X \to Y$ be a morphism of schemes. Assume that~$f$ is proper, flat, and of finite presentation. Then the set of $y \in Y$ such that the fiber $X_y \coloneqq f^{-1}(y)$ is geometrically irreducible is Zariski open in~$Y$.
\end{theorem}

\begin{proof}
See~\cite[Th\'eor\`eme 12.2.4(viii)]{EGAIV.3}.
\end{proof}

\begin{remark}
\label{rem:reduction}
In the event that~$Y$ is a locally Noetherian scheme, one can replace the condition ``of finite presentation'' with ``of finite type''\cite[\href{https://stacks.math.columbia.edu/tag/01TX}{Tag 01TX}]{stacks-project}. However, this condition is in turn subsumed by the properness condition (by definition of properness!).
\end{remark}

Since the coordinate ring~$A$ of the affine space $\Aff_\Q^m$ is Noetherian, the scheme $\Aff_\Q^m$ is locally Noetherian.  By Remark~\ref{rem:reduction}, to deduce Theorem~\ref{thm:precisegeomirreduc} from Theorem~\ref{thm:GenericIrreducibilityEGA}, we must show that the slowness surface fibration $f \colon \bS \to \Aff_\Q^{m}$ is a proper, flat morphism.  We say a few words about what these conditions mean first.

In algebraic geometry, the notion of properness mimics the analogous notion between complex analytic spaces: the preimage of a compact set is compact. In particular, a proper morphism takes closed sets to closed sets; see~\cite[\href{https://stacks.math.columbia.edu/tag/01W0}{Section 01W0}]{stacks-project}. Flatness is an algebraic condition that, in conjunction with properness and local Noetherianity of the target, guarantees the nonempty fibers of~$f$ vary nicely (e.g., they all have the same Euler characteristic); see \cite[\href{https://stacks.math.columbia.edu/tag/01U2}{Section 01U2}]{stacks-project}. To prove that the slowness surface fibration $f\colon \bS \to \Aff_\Q^m$ is flat, we shall use the ``miracle flatness'' criterion.

\begin{theorem}[Miracle Flatness]
\label{thm:miracleflatness}
Let $f\colon X \to Y$ be a morphism of finite type, equidimensional schemes over a field. Suppose that~$X$ is Cohen-Macaulay,~$Y$ is regular, and the fibers of~$f$ have dimension $\dim X - \dim Y$.  Then~$f$ is a flat morphism.
\end{theorem}

\begin{proof}
See~\cite[26.2.11]{Vakil}.
\end{proof}

\begin{proposition} 
\label{prop:HypothesisCheck}
The slowness surface fibration $f\colon \bS \to \Aff_\Q^{m}$ is a proper flat morphism. If $K/\Q$ is a field extension, the same conclusion holds for the base-extension $f_K\colon S_K \to \Aff_K^m$.
\end{proposition}

\begin{proof}
First, we prove that~$f$ is proper.  The scheme $\PP^n_\Q$ being projective, its structure morphism $\PP^n_\Q \to \Spec \Q$ is proper.  Consider the fibered product diagram
\[
\xymatrix{
\Aff_\Q^{m} \times_{\Spec \Q} \PP_\Q^n \ar[r]^{\pi_2} \ar[d]^{\pi_1} & \PP_\Q^n \ar[d] \\
\Aff_\Q^{m} \ar[r] & \Spec \Q
}
\]
Proper morphisms are stable under base change~\cite[II Corollary 4.8(c)]{Hartshorne}, and hence~$\pi_1$ is proper. Closed immersions being proper~\cite[II Corollary 4.8(a)]{Hartshorne}, the morphism $\iota\colon S \into \PP^n_A$ is also proper.  Finally, a composition of proper morphisms is proper~\cite[II Corollary 4.8(b)]{Hartshorne}, whence $f = \pi_1\circ\iota$ is proper.  

Next, we show that the morphism $f\colon \bS \to \Aff_\Q^{m}$ is flat via Theorem~\ref{thm:miracleflatness}.  The schemes~$\bS$ and $\Aff_\Q^{m}$ are of finite type over a field and equidimensional ($\bS$ is a hypersurface in~$\PP^n_A$), so it suffices to verify that~$\bS$ is a Cohen--Macaulay scheme, that $\Aff_\Q^{m}$ is regular, and that the fibers of~$f$ all have dimension $n-1$.  The surface~$\bS$ is a hypersurface in a projective space, so it is a local complete intersection, hence a Cohen-Macaulay scheme~\cite[\href{https://stacks.math.columbia.edu/tag/00SA}{Tag 00SA}]{stacks-project}. The affine space $\Aff_\Q^{m}$ is smooth, hence regular; the fibers of~$f$ are all hypersurfaces of $\PP^n_{\Q}$, because the coefficient of~$p_0^{2n}$ in the defining equation of every fiber is $\pm 1$, and hence all have dimension $n-1$. Hence, the morphism~$f$ is flat.

The claim for the base-extension $f_K\colon S_K \to \Aff_K^m$ follows either by replacing~$\Q$ with~$K$ in the above arguments, or by noting that proper and flat are properties of morphisms that are stable under base-extension (see, e.g.,~\cite[\href{https://stacks.math.columbia.edu/tag/01U9}{Lemma 01U9}]{stacks-project}).
\end{proof}

\begin{proof}[Proof of Theorem~\ref{thm:precisegeomirreduc}]
The conclusion that $\Irr(f)$ is Zariski open in $\Aff_\R^m$ follows from Theorem~\ref{thm:GenericIrreducibilityEGA} and~Proposition~\ref{prop:HypothesisCheck}, taking into account Remark~\ref{rem:reduction}. It follows that $\Irr(f)$ is either empty, or it is the complement of a proper closed subset of $\Aff_\R^m$. Such a set is determined by an ideal $I \subseteq \R[b_{ij} : 1 \leq i,j \leq \pairs(n) ]$~\cite[Corollary~II.5.10]{Hartshorne}. The base-extension $I_\C = I\otimes_\R \C$ determines a proper closed subset of $\Aff_\C^m$, whose finitely many irreducible components have dimension $\leq m - 1$. This closed subset descends to finitely many irreducible components in $\Aff_\R^m$, consisting of complex conjugate pairs of irreducible varieties in $\Aff_\C^m$.
\end{proof}

%%%%%%%%%%%
\subsection{Ex uno plura}
All our work so far does not preclude the possibility that the Zariski open subset $\Irr(f)$ of $\Aff_\R^m$ defined in Theorem~\ref{thm:precisegeomirreduc} is empty! We verify that this is not in the context of Examples~\ref{ex:n=2}--\ref{ex:n=3ortho} by exhibiting explicit slowness polynomials of the appropriate form that are irreducible over~$\C$.  We use a standard arithmetic trick: reduction modulo a prime. The principle involved is simple: if a polynomial $F(x_0,\dots,x_n)$ with coefficients in~$\Z$ factors nontrivially, then it also factors when we reduce its coefficients modulo any prime~$p$. Thus, if a polynomial with coefficients in~$\Z$ is irreducible when considered over the finite field~$\F_p$, then it must be irreducible over~$\Z$. This principle is extraordinarily useful, because by finiteness of~$\F_p$ checking whether the reduction $\bar F(x_0,\dots,x_n) \in \F_p[x_0,\dots,x_n]$ is irreducible is a finite, fast computation. Guaranteeing that the polynomial remains irreducible when considered over~$\C$ requires working over a finite extension~$\F_{p^d}$ of~$\F_p$ with controlled degree~$d$. We make this idea explicit in the following lemma, whose proof we include for lack of a good reference.

\begin{lemma}
\label{lem:IrreducibilityCriterion}
Let $F(x_0,\dots,x_n) \in \Z[x_0,\dots,x_n]$ be a homogeneous polynomial of degree~$d$. Suppose there is a prime~$p$ such that the reduction $\bar F(x_0,\dots,x_n) \in \F_p[x_0,\dots,x_n]$ of~$F$ modulo~$p$ is irreducible in the finite field~$\F_{p^d}$ of cardinality~$p^{d}$. Then $F(x_0,\dots,x_n)$ is irreducible in $\C[x_0,\dots,x_n]$.
\end{lemma}

\begin{proof}
By~\cite[\href{https://stacks.math.columbia.edu/tag/020J}{Tag 020J}]{stacks-project}, to prove that the polynomial~$F$ is irreducible over~$\C$, it suffices to show that it is irreducible in $\overline{\Q}[x_0,\dots,x_n]$, where $\overline{\Q}$ denotes a fixed algebraic closure of~$\Q$. The field~$\overline{\Q}$ consists of all algebraic numbers: the roots of single-variable polynomials with rational coefficients; it is countable.

There is a Galois field extension
$K/\Q$ of \emph{finite} degree where~$F$ already factors into $\overline{\Q}$-irreducible polynomials. To see this, note that each coefficient of each factor in a~$\overline{\Q}$ factorization is an algebraic number, hence has finite degree over~$\Q$; we can let~$K$ be the Galois closure of the field obtained from~$\Q$ by adjoining all the coefficients of all the factors of~$F$ over~$\overline{\Q}$ (see also~\cite[\href{https://stacks.math.columbia.edu/tag/04KZ}{Tag 04KZ}]{stacks-project}).

Let~$\frakp$ be a prime ideal in the ring of integers~$\calO_K$ of~$K$ lying over~$p$, i.e., $\frakp\cap \Z = (p)$.  The field $\F_\frakp \coloneqq \calO_K/\frakp\calO_K$ is a finite field extension of~$\F_p$. Let $\bar F = g_1,\dots,g_m$ be a factorization of~$\bar F$ in $\F_\frakp[x_0,\dots,x_n]$.  The Galois group $G \coloneqq \Gal(\F_\frakp/\F_p)$ acts on the set $\{g_1,\dots,g_m\}$. The orbits of this action correspond to the irreducible factors of the reduction of $\bar F \in \F_p[x_0,\dots,x_n]$. This reduction is irreducible, because by hypothesis~$\bar F$ is irreducible over the larger field~$\F_{p^d}$, so the action of~$G$ on $\{g_1,\dots,g_m\}$ is transitive. It follows that the factors $\{g_1,\dots,g_m\}$ of~$\bar F$ must all have the same degree, and hence $m \mid d$. By the orbit-stabilizer theorem, the stabilizer $H_{g_i} \leq G$ of~$g_i$ has index~$m$; it is a normal subgroup of~$G$ because~$G$ is cyclic. By Galois theory, the polynomial~$\bar F$ already factors over the fixed field $K(\frakp)^{H_{g_i}} \isom \F_{p^{m}}$ as $g_i\cdot h_i$ for some $h_i \in \F_{p^{m}}[x_0,\dots,x_n]$. However, by hypothesis, the polynomial $\bar F \in \F_p[x_0,\dots,x_n]$ is irreducible over~$\F_{p^{d}}$, and hence is irreducible over~$\F_{p^{m}}$, because $\F_{p^m} \subset \F_{p^d}$ as $m \mid d$. This implies that $m = 1$, i.e.,~$\bar F$ is irreducible in $\F_\frakp[x_0,\dots,x_n]$, and therefore~$F$ irreducible in $K[x_0,\dots,x_n]$. By definition of the field~$K$, we conclude that~$F$ is irreducible in $\overline{\Q}[x_0,\dots,x_n]$.
\end{proof}

\begin{remark}
The hypothesis that $F \in \Z[x_0,\dots,x_n]$ is homogeneous can be weakened. We used this hypothesis tacitly above: we assumed that the reduction $\bar F \in \F_\frakp[x_0,\dots,x_n]$ has degree~$d$. This is certainly the case if~$F$ is homogeneous and~$\bar F$ is irreducible (and hence nonzero).
\end{remark}

\begin{example}
\label{prop:2D-example}
Let $n=2$. Using the notation of Example~\ref{ex:n=2}, consider the stiffness tensor with components
\[
b_{11} = 20, \quad
b_{12} = 39, \quad
b_{13} = -65, \quad
b_{22} = -16, \quad
b_{23} = -87, \quad
b_{33} = 30.
\]
The corresponding homogenized slowness polynomial
\begin{equation*}
\label{eq:2Dexample}
\tilde P(\vect{p}) = 
-3625p_1^4 + 1590p_1^3p_2 + 7129p_1^2p_2^2 - 50p_1^2p_0^2 + 8866p_1p_2^3
+ 304p_1p_2p_0^2 - 8049p_2^4 - 14p_2^2p_0^2 + p_0^4
\end{equation*}
is irreducible over~$\C$: apply Lemma~\ref{lem:IrreducibilityCriterion} with $d = 4$ and $p = 7$: a {\tt magma} calculation shows that the reduction of this polynomial modulo~$7$ is irreducible in the finite field~$\F_{7^4}$; see~\cite{this-paper}.
\end{example}

\begin{example}
    \label{ex:olivineirreduc}
    Let $n=3$. Consider the orthorhombic stiffness tensor obtained by rounding out values for the stiffness tensor of olivine~\cite{olivine}, a common mineral in the Earth's mantle:
    \begin{equation}
    \label{eq:orthorhombicparameterexample}
        \begin{split}
                b_{11} = 321, \quad b_{12} &= 68, \quad b_{13} = 72, \quad b_{22} = 197, \quad b_{23} = 77, \\ b_{33} &= 234, \quad b_{44} = 64, \quad b_{55} = 77, \quad b_{66} = 79.
        \end{split}
    \end{equation}
    Its corresponding homogenized slowness polynomial is
    \begin{equation}
    \label{eq:orthorhombicslownesspolyexample}
        \begin{split}
            \tilde P(\vect{p}) &= 
        1952643p_1^6 + 5308889p_1^4p_2^2 + 6230406p_1^4p_3^2 - 56159p_1^4p_0^2 + 4261967p_1^2p_2^4 \\
        &\quad+ 9884047p_1^2p_2^2p_3^2 - 
    94721p_1^2p_2^2p_0^2 + 5189310p_1^2p_3^4 - 108883p_1^2p_3^2p_0^2 + 477p_1^2p_0^4 \\
    &\quad+ 996032p_2^6 + 3365543p_2^4p_3^2 -
    33227p_2^4p_0^2 + 3517205p_2^2p_3^4 - 73952p_2^2p_3^2p_0^2 \\
    &\quad+ 340p_2^2p_0^4 + 1153152p_3^6 - 37922p_3^4p_0^2 + 
    375p_3^2p_0^4 - p_0^6.
        \end{split}
    \end{equation}
    which is irreducible over~$\C$ by Lemma~\ref{lem:IrreducibilityCriterion}, applied with $d = 6$ and $p = 5$; see~\cite{this-paper}.
\end{example}

\begin{proof}[Proof of Theorem~\ref{thm:GenericIrreducibility}]
By Theorem~\ref{thm:precisegeomirreduc} we know that the subset $\Irr(f)$ of the parameter space of stiffness tensors whose corresponding homogenized slowness polynomials are irreducible over~$\C$ is a Zariski open subset of $\Aff_\R^m$. Since $\Aff_\R^m$ is irreducible, a Zariski open subset is dense, as long as it is not empty.  Example~\ref{prop:2D-example} shows that $\Irr(f)$ is not empty when $n = 2$; Example~\ref{ex:olivineirreduc} shows that $\Irr(f)$ is not empty for slowness surface fibration corresponding to orthorhombic stiffness tensors. Since $\orthorhombic \subset \monoclinic \subset \stiffnesstensor{3}$, Example~\ref{ex:olivineirreduc} also furnishes the required non-emptyness statement in the cases of monoclinic stiffness tensors and fully anisotropic stiffness tensors (i.e., the case $n = 3$).
\end{proof}

%%%%%%%%%%%
\subsection{Generic reconstruction of stiffness tensors}
\label{sec:recoverability}
In this section, we prove Theorem~\ref{thm:polynomial-to-tensor}.  While the \emph{proof} uses elaborate ideas from algebraic geometry, we may perform the reconstruction of a stiffness tensor from a particular slowness polynomial \emph{quickly in practice}, using simple ideas from the theory of Gr\"obner bases (see~\S\ref{ss:2DGB} and~\cite{this-paper}).

We begin with the two-dimensional case. In our first proof of Theorem~\ref{thm:polynomial-to-tensor}(\ref{item:2D}), we use a criterion of Eisenbud and Ulrich~\cite{EisenbudUlrich} to determine when a rational map $g\colon \PP^r \dasharrow \PP^s$ is birational onto its image.  Write $R = \Q[x_1,\dots,x_{r+1}]$ for the coordinate ring of $\PP^r$.  The map $g$ is given by a sequence of homogeneous polynomials $g_{1},\dots,g_{s+1}$ in $R$ of the same degree (assumed, without loss of generality, to have no global common factor).  We let $I = \langle g_1,\dots,g_{s+1}\rangle \subset R$ be the ideal generated by these homogeneous polynomials.  A point $y = (c_1:\,\cdots\,:c_{s+1}) \in \PP^s(\Q)$ gives rise to the ideal $I_y\subset R$ generated by the $2\times 2$ minors of the matrix
\[
\begin{pmatrix}
g_1(x_1,\dots,x_{r+1}) & \cdots & g_{s+1}(x_{1},\dots,x_{r+1}) \\
c_1 & \cdots & c_{s+1}
\end{pmatrix}
\]
The ideal quotient $I_y : I$ is called the \defi{generalized row ideal} corresponding to $y$.  An ideal $J \subset R$ has \defi{linear type}, or is \defi{linear} if the natural map
\[
\Sym(J) \to \mathcal{R}(J)
\]
from the symmetric algebra of $J$ to the Rees algebra of $J$ is an isomorphism.

\begin{proposition}[{\cite[Proposition~3.1(c)]{EisenbudUlrich}}]
\label{prop:EisenbudUlrich}
If there is a point $y \in \PP^s(\Q)$ such that the ideal $I_y : I$ is linear of codimension $r$, then the map $g$ is birational onto its image.
\qed
\end{proposition}

\begin{proof}[First Proof of Theorem~\ref{thm:polynomial-to-tensor}(\ref{item:2D})]
Using the notation of~\S\ref{sec:alg-principles}, we define the rational map~\eqref{eq:ratmap2D} between complex projective spaces
\begin{equation*}
\begin{split}
g\colon \PP^6_{(b_{11}:\,\cdots\,:b_{33}:r)} &\dasharrow \PP^7_{(\tilde c_1:\,\cdots\,:\tilde c_8)}, \\
(b_{11}:\,\cdots\,:b_{33}:r) &\mapsto \left(b_{11}b_{33} - b_{13}^2:\,\cdots\,:-r(b_{22} + b_{33})\right).
\end{split}
\end{equation*}
Let 
\begin{equation}
\label{eq:Itoblowup}
I = \langle b_{11}b_{33} - b_{13}^2,\dots,-r(b_{22} + b_{33}) \rangle \subset \Q[b_{11},\dots,b_{33},r],
\end{equation}
and let 
\[
y = (-3625:1590:7129:-50:8866:304:-8049:-14) \in \PP^7(\Q).
\] 
Then we can verify in Macaulay2 that the ideal quotient $I_y\colon I$ is linear of codimension 6. By Proposition~\ref{prop:EisenbudUlrich}, the map $g$ is birational onto its image.  A computation verifying this claim is included in~\cite{this-paper}. By our discussion in \S\ref{sec:uniquereconstruction}, this is enough to show the map $\spolymap\colon\stiffnesstensor{2} \to \R[p_1,p_2]$ is generically injective.
\end{proof}

While in principle, one could use Proposition~\ref{prop:EisenbudUlrich} to prove the analog of Theorem~\ref{thm:polynomial-to-tensor}(\ref{item:2D}) in dimension~$n=3$, the required symbolic computation seems to be out of current reach.  Nevertheless, we included the proof above in dimension $n=2$ in the hopes it might be useful in other inverse-problem situations.

\begin{proof}[Second Proof of Theorem~\ref{thm:polynomial-to-tensor}(\ref{item:2D})]
We use the same rational map $g\colon \PP^6 \dasharrow \PP^7$ as in the first proof of the theorem. The closed subset~$\Pi\subset \PP^6$ where~$g$ is not defined is $2$-dimensional, although we will not use this fact explicitly.  Let $X\coloneqq \Bl_\Pi(\PP^6)$ be the blow-up of~$\PP^6$ along~$\Pi$ \cite[Example~II.7.17.3]{Hartshorne}. This scheme comes with a morphism $\pi \colon X \to \PP^6$ such that the composition 
\[
h \coloneqq g\circ\pi\colon X \to \PP^7
\]
is a proper morphism, and such that $X\setminus \pi^{-1}(\Pi) \isom \PP^6 \setminus \Pi$. Let~$Y$ be the image of~$h$, so that the map $h\colon X \to Y$ is a surjective proper morphism. Properness ensures that the \defi{fiber dimension function}
\begin{align*}
d\colon Y &\to \R \\
y &\mapsto \dim h^{-1}(y)
\end{align*}
is upper semi-continuous, i.e., for each $x\in \R$ the set $d^{-1}((-\infty,x))$ is Zariski open~\cite[\href{https://stacks.math.columbia.edu/tag/0D4I}{Tag 0D4I}]{stacks-project}.  In particular, if there is a point $y \in Y$ such that $\dim h^{-1}(y) = 0$, then there is a nonempty Zariski open subset $V\subseteq Y$ over which all fibers are $0$-dimensional.  (Note that both~$X$ and~$Y$ are irreducible varieties.) Assuming such a point $y \in Y$ exists, the fibers over~$V$ have finite cardinality, so the induced morphism $h\colon U \coloneqq h^{-1}(V) \to V$ is quasi-finite.  It is also proper, as it is a base-extension of a proper morphism.  A proper, quasifinite morphism is finite~\cite[\href{https://stacks.math.columbia.edu/tag/02OG}{Tag 02OG}]{stacks-project}.  Finally, the fiber degree function is also an upper semi-continuous function on the target of a finite morphism~\cite[14.3.4]{Vakil}. Summing up, if there is a point $u \in U$ such that $h^{-1}(h(u))$ consists of a \emph{single} point, then upper semi-continuity of degree implies there is a Zariski open subset $V' \subseteq V$ such that, for all $y \in V'$, the fiber $h^{-1}(y)$ consists of exactly one point, proving generic injectivity for the morphism $h$ onto its image.

We use the point $u \coloneqq \pi^{-1}(v)$, where
\[
v = (b_{11}:b_{12}:b_{13}:b_{22}:b_{23}:b_{33}:r) = (20:39:-65:-16:-87:30:1) \in \PP^6(\Q)
\]
Since $v \in \PP^6\setminus \Pi$, we know that $\pi^{-1}(v)$ consists of exactly one point. Its image in $\PP^7$ is
\[
y \coloneqq h(u) = g\circ \pi(u) = g(v) = (3625:-1590:-7129:50:-8866:-304:8049:14),
\]
and we want to show that $h^{-1}(y)$ consists of exactly one (closed) point, scheme-theoretically. A computation verifying this claim is included in~\cite{this-paper}, but we explain the idea here for the benefit of readers coming at this paper from an inverse problems background.

\medskip

Let $I \subseteq \Q[b_{11},\dots,b_{33},r] =: R$ be the ideal~\eqref{eq:Itoblowup} in the first proof of the theorem. The blow-up of $\PP^6$ along the vanishing of $I$ is $\Proj \mathcal{R}(I)$, where $\mathcal{R}(I)$ is the Rees algebra of $R$, which is
\[
\mathcal{R}(I) = R[tI] = \oplus_{i = 0}^\infty t^iI^i = R \oplus It \oplus I^2t^2 \oplus t^3I^3 \oplus \cdots
\]
see, for example~\cite[\S~IV.2]{EisenbudHarris}. Note that $I^0 = R$, and we should view $\mathcal{R}(I)$ as a bi-graded $\Q$-algebra, where the variables $b_{ij}$ and~$r$ have degree $(1,0)$, and the variable $t$ has degree $(0,1)$. So, for example, elements in $I^nt^n$ have degree $(n,n)$.
A {\tt magma} calculation shows that 
\[
\mathcal{R}(I) \isom R[\tilde c_1,\dots,\tilde c_8]/J \isom \Q[b_{11},\dots,b_{33},r,\tilde c_1,\dots,\tilde c_8]/J,
\]
where
\begingroup
\allowdisplaybreaks
\begin{align*}
    J = \langle&2b_{13}\tilde c_8 + 2b_{23}\tilde c_8 - b_{22}\tilde c_5 - b_{33}\tilde c_5, \\
    &b_{11}\tilde c_8 + b_{33}\tilde c_8 - b_{22}\tilde c_4 - b_{33}\tilde c_4, \\
    &2b_{33}\tilde c_8 + b_{13}\tilde c_5 - b_{23}\tilde c_5 - 2b_{33}\tilde c_4 + 2r\tilde c_7 - 2r\tilde c_1, \\
    &2b_{23}\tilde c_8 + b_{12}\tilde c_5 - b_{22}\tilde c_5 - 2b_{23}\tilde c_4 - r\tilde c_5 - r\tilde c_2, \\
    &b_{11}\tilde c_5 + b_{33}\tilde c_5 - 2b_{13}\tilde c_4 - 2b_{23}\tilde c_4, \\
    &2b_{11}\tilde c_7 - b_{13}\tilde c_5 - 2b_{22}\tilde c_1 + b_{23}\tilde c_2, \\
    &-4b_{13}\tilde c_7 + b_{12}\tilde c_5 + 2b_{33}\tilde c_5 + b_{22}\tilde c_2 - 2b_{23}\tilde c_3, \\
    &b_{11}\tilde c_5 - 4b_{23}\tilde c_1 + b_{12}\tilde c_2 + 2b_{33}\tilde c_2 - 2b_{13}\tilde c_3, \\
    &b_{23}^2\tilde c_8 - b_{22}b_{33}\tilde c_8 - b_{22}r\tilde c_7 - b_{33}r\tilde c_7, \\
    &-2b_{13}b_{33}\tilde c_8 - 2b_{23}b_{33}\tilde c_8 + b_{23}^2\tilde c_5 + b_{33}^2\tilde c_5 - 2b_{13}r\tilde c_7 - 2b_{23}r\tilde c_7, \\
    &-b_{11}b_{33}\tilde c_8 - b_{33}^2\tilde c_8 + b_{23}^2\tilde c_4 + b_{33}^2\tilde c_4 - b_{11}r\tilde c_7 - b_{33}r\tilde c_7, \\
    &-2b_{12}b_{23}\tilde c_8 - 2b_{22}b_{23}\tilde c_8 + b_{22}^2\tilde c_5 + b_{22}b_{33}\tilde c_5 + b_{22}r\tilde c_5 + b_{33}r\tilde c_5, \\
    &-2b_{11}b_{23}\tilde c_8 - 2b_{23}b_{33}\tilde c_8 + b_{11}b_{22}\tilde c_5 + b_{22}b_{33}\tilde c_5 - 2b_{12}b_{23}\tilde c_4 + 2b_{23}b_{33}\tilde c_4 + b_{11}r\tilde c_5 + b_{33}r\tilde c_5, \\
    &b_{12}^2\tilde c_8 - b_{11}b_{22}\tilde c_8 - 2b_{13}b_{23}\tilde c_8 + 2b_{12}b_{33}\tilde c_8 - b_{22}r\tilde c_3 - b_{33}r\tilde c_3, \\
    &b_{11}^2\tilde c_8 + b_{11}b_{33}\tilde c_8 - 2b_{12}b_{33}\tilde c_8 - 2b_{22}b_{33}\tilde c_8 + b_{22}b_{23}\tilde c_5 + b_{23}b_{33}\tilde c_5 - b_{12}^2\tilde c_4 - b_{11}b_{33}\tilde c_4 - 2b_{12}r\tilde c_7 \\
    &\quad\quad - 2b_{22}r\tilde c_7 - 2b_{33}r\tilde c_7 + b_{23}r\tilde c_5 + 2b_{12}r\tilde c_1 - b_{13}r\tilde c_2 + b_{11}r\tilde c_3 + b_{33}r\tilde c_3, \\
    &2b_{13}b_{22}\tilde c_7 - 2b_{12}b_{23}\tilde c_7 + b_{23}^2\tilde c_5 - b_{22}b_{33}\tilde c_5, \\
    &2b_{13}^2\tilde c_7 - b_{13}b_{33}\tilde c_5 - 2b_{23}^2\tilde c_1 + b_{23}b_{33}\tilde c_2, \\
    &2b_{12}b_{13}\tilde c_7 - 2b_{11}b_{23}\tilde c_7 + 4b_{13}b_{33}\tilde c_7 - b_{12}b_{33}\tilde c_5 - 2b_{33}^2\tilde c_5 - b_{23}^2\tilde c_2 + 2b_{23}b_{33}\tilde c_3, \\
    &b_{12}^2\tilde c_7 - b_{11}b_{22}\tilde c_7 - 2b_{13}b_{23}\tilde c_7 + 2b_{12}b_{33}\tilde c_7 - b_{23}^2\tilde c_3 + b_{22}b_{33}\tilde c_3, \\
    &2b_{11}b_{13}\tilde c_7 - b_{11}b_{33}\tilde c_5 - 2b_{12}b_{23}\tilde c_1 + b_{13}b_{23}\tilde c_2, \\
    &2b_{12}b_{13}\tilde c_1 - 2b_{11}b_{23}\tilde c_1 - b_{13}^2\tilde c_2 + b_{11}b_{33}\tilde c_2, \\
    &2b_{11}^2\tilde c_7 - 2b_{12}^2\tilde c_1 - 4b_{12}b_{33}\tilde c_1 + b_{12}b_{13}\tilde c_2 + b_{11}b_{23}\tilde c_2 +  2b_{13}b_{33}\tilde c_2 - 2b_{11}b_{33}\tilde c_3\rangle
\end{align*}
\endgroup
Under this isomorphism, we should think of $\mathcal{R}(I)$ as a bi-graded $\Q$-algebra where the variables $b_{ij}$ and~$r$ have degree $(1,0)$, and the variables $\tilde c_i$ have degree $(0,1)$. This allows us to think of $X = \Proj \mathcal{R}(I)$ as a subscheme of the product $\PP^6 \times \PP^7$, and the maps $\pi \colon X \to \PP^6$ and $h\colon X \to \PP^7$ are the projection maps onto each of the factors, induced by the applying the $\Proj$ construction to the natural graded ring maps $\Q[b_{11},\dots,b_{33},r] \to \mathcal{R}(I)$ and $\Q[\tilde c_1,\dots,\tilde c_8] \to \mathcal{R}(I)$.

The fiber $h^{-1}(y)$ can now be computed as the subscheme $V(J + J') \subset \PP^6\times \PP^7$, where $J'$ is the ideal in $\Q[b_{11},\dots,b_{33},r,\tilde c_1,\dots,\tilde c_8]$ generated by the $2\times 2$ minors of the matrix
\[
\begin{pmatrix}
    3625 & -1590 & \cdots & 14 \\
    \tilde c_1 & \tilde c_2 & \cdots & \tilde c_8
\end{pmatrix}
\]
A {\tt magma} calculation reveals
\begin{align*}
J + J' = \langle &b_{11} - 20r, b_{12} - 39r, b_{13} + 65r,b_{22} + 16r, b_{23} + 87r,b_{33} - 30r \\
&14\tilde c_1 - 3625\tilde c_8,14\tilde c_2 + 1590\tilde c_8,14\tilde c_3 + 7129\tilde c_8,14\tilde c_4 - 50\tilde c_8, \\
&14\tilde c_5 + 8866\tilde c_8,14\tilde c_6 + 304\tilde c_8,14\tilde c_7 - 8049\tilde c_8\rangle.
\end{align*}
This calculation confirms that the scheme $h^{-1}(y) \subset X$ consists of exactly one point in $\PP^6\times \PP^7$, with coordinates
\[
\left( (20:39:-65:-16:-87:30:1),(3625:-1590:-7129:50:-8866:-304:8049:14)\right).
\]

We conclude that the map $h\colon X \to Y$ is generically injective.  Note that the locus where $r = 1$ is the distinguished dense open affine chart $\Aff^6 \isom \{r \neq 0\} = D_+(r) \subset \PP^6$, and that~$h$ and~$g$ coincide on $D_+(r)\cap (\PP^6 \setminus \Pi)$, so~$h$ is still generically injective after ``dehomogenizing~$r$''. \end{proof}

The second proof of Theorem~\ref{thm:polynomial-to-tensor}(\ref{item:2D}) readily generalizes to parts~(\ref{item:3Dortho}) and~(\ref{item:3Dmono}) of the theorem in a way that the first proof cannot: the method of the first proof can only show generic injectivity of a map, whereas the remaining two cases of the theorem requiring proving $\spolymap|_{\orthorhombic}$ and $\spolymap|_{\monoclinic}$ are generically $4$-to-$1$ and $2$-to-$1$, respectively.

\begin{proof}[Proof of Theorem~\ref{thm:polynomial-to-tensor}(\ref{item:3Dortho})]
Recall an orthorhombic stiffness tensor has $9$ parameters $b_{11}$, $b_{12}$, $b_{13}$, $b_{22}$, $b_{23}$, $b_{33}$, $b_{44}$, $b_{55}$, and $b_{66}$, and that its corresponding slowness polynomial~\eqref{eq:genericorthorhombic} has $20$ coefficients $c_1,\dots,c_{20}$ (although $c_{20} = -1$ always, so we ignore it in our discussion). 
Homogenizing these coefficients with respect to an auxiliary variable $r$ as we did in~\eqref{eq:2Dreconstruction}, we arrive at a rational map
\[
\PP^9_{(b_{11}:\,\cdots\,:b_{66},r)} \dasharrow \PP^{18}_{(\tilde c_1:\,\cdots\,:\tilde c_{19})}
\]
where we could apply reasoning similar to the second proof of Theorem~\ref{thm:polynomial-to-tensor}(\ref{item:2D}). However, a small amount of linear algebra simplifies the setup, reducing the computational cost required to run the argument. We must blow up an ideal $I \subset \Q[b_{11},\dots,b_{66},r]$, generated by the expressions $\tilde c_1,\dots,\tilde c_{19}$ in terms of the variables $b_{11},\dots,b_{66},r$.  However, three of the original coefficients $c_1,\dots,c_{19}$ are quite special:
\begin{equation}
\label{eq:linearrelations}
\begin{split}
    c_{10} &= b_{11} + b_{55} + b_{66}, \\
    c_{16} &= b_{22} + b_{44} + b_{66}, \\
    c_{19} &= b_{33} + b_{44} + b_{55}.
\end{split}    
\end{equation}
These linear relations imply that for a general set of stiffness parameters, the quantities $b_{33}$, $b_{44}$, and $b_{66}$ are determined uniquely by the values of $b_{11}$, $b_{22}$, $b_{55}$, $c_{10}$, $c_{16}$ and $c_{19}$. We can use the generators $b_{11} + b_{55} + b_{66}$, $b_{22} + b_{44} + b_{66}$, and $b_{33} + b_{44} + b_{55}$ of $I$ to simplify the rational map $\PP^9 \dasharrow \PP^{18}$ by eliminating $b_{33}$, $b_{44}$, and $b_{66}$ and looking at the map
\begin{align*}
g \colon\PP^6 &\dasharrow \PP^{15} \\
(b_{11}:b_{12}:b_{13}:b_{22}:b_{23}:b_{55}:r) &\mapsto \left(\tilde c_1:\,\cdots\,:\widehat{\tilde c}_{10}:\,\cdots\,:\widehat{\tilde c}_{16}:\,\cdots\,:\tilde c_{18},\widehat{\tilde c}_{19}\right),
\end{align*}
where, for example, the notation $\widehat{\tilde c}_{10}$ indicates that the coordinate $\tilde c_{10}$ has been omitted.

Now we proceed as in the proof of Theorem~\ref{thm:polynomial-to-tensor}(\ref{item:2D}): after resolving the indeterminacy locus\footnote{This locus has dimension~$1$, as one can verify with {\tt magma}, for example.} $\Pi$ of~$g$ through a blow-up process to get a surjective proper morphism $h\colon X \to Y$, upper semi-continuity of fiber dimension together with upper semi-continuity of degree for finite morphisms show there is a Zariski open subset of~$Y$ over which all fibers consist of \emph{at most four points}: consider the point
\[
v = (b_{11}:b_{12}:b_{13}:b_{22}:b_{23}:b_{55}:r) =
(321:68:72:197:77:77:1) \in \PP^{6}(\Q),
\]
coming from the stiffness parameters for olivine. It lies in $\PP^{6}\setminus\Pi$, so its preimage under the blow-up map $\pi\colon X \to \PP^{6}$ consists of a single point $u$ such that
\begin{align*}
y :&= h(u) = g(v) \\
&= ( -9837366: -16998301: 46961097: 133687: 16623679: -86045237: -81787:-12687922:\\
&\quad -54486: -15759606: 39525271: 118807: 15182471: 30357: -4302606: 61807 ) \in \PP^{15}(\Q).
\end{align*}
A Gr\"obner basis calculation included in~\cite{this-paper} shows that $h^{-1}(y)\subset \PP^6(\Q)$ consists of four closed points, namely
\begin{align*}
(321:68:72:197:77:77:1),&\quad (321:728:-226:197:77:77:1),\\
(321:68:-226:197:-479:77:1),&\quad (321:728:72:197:-479:77:1),
\end{align*}
which includes the original stiffness parameters and three other anomalous companions.

To complete the proof, we show that a non-empty fiber of the map $g$ has \emph{at least} four points. The above calculation suggests that anomalous companions can differ in the stiffnesses $b_{12}$, $b_{13}$, and $b_{23}$. Inspecting the coefficients of~\eqref{eq:genericorthorhombic} we see that
\begin{align*}
c_7 &= b_{12}^2 + 2b_{12}b_{66} -b_{11}b_{22} - b_{11}b_{44}  - b_{22}b_{55} - b_{44}b_{66} - b_{55}b_{66}, \\
\label{eq:c8}
c_9 &= b_{13}^2 + 2b_{13}b_{55} -b_{11}b_{33} - b_{11}b_{44} - b_{33}b_{66} - b_{44}b_{55} - b_{55}b_{66}, \\
c_{15} &= b_{23}^2 + 2b_{23}b_{44} - b_{22}b_{33} - b_{22}b_{55} - b_{33}b_{66} - b_{44}b_{55} - b_{44}b_{66}.
\end{align*}
Considering these equations as quadratic in $b_{12}$, $b_{13}$, and $b_{23}$, respectively, suggests that each of these stiffnesses can take on two possible values $b_{ij}$ and $b_{ij}^*$, related by the equations
\begin{align*}
b_{12} + b_{12}^* &= -2b_{66}, \\ 
b_{13} + b_{13}^* &= -2b_{55}, \\ 
b_{23} + b_{23}^* &= -2b_{44}.
\end{align*}
These equations are coupled by the linear relations~\eqref{eq:linearrelations}, so only four possible triples can occur as stiffness parameters:
\begin{equation}
\label{eq:anomalouscompanions}
(b_{12},b_{13},b_{23}),\quad
(b_{12},b_{13}^*,b_{23}^*),\quad
(b_{12}^*,b_{13}^*,b_{23}),\quad
(b_{12}^*,b_{13},b_{23}^*),    
\end{equation}
Now we can simply verify that the four points
\begin{align*}
(b_{11}:b_{12}:b_{13}:b_{22}:b_{23}:b_{55}:r),&\ \ (b_{11}:-2b_{66}-b_{12}:-2b_{55}-b_{13}:b_{22}:b_{23}:b_{55}:r),\\
(b_{11}:b_{12}:-2b_{55}-b_{13}:b_{22}:-2b_{44}-b_{23}:b_{55}:r),&\ \  (b_{11}:-2b_{66}-b_{12}:b_{13}:b_{22}:-2b_{44}-b_{23}:b_{55}:r),
\end{align*}
have the same image under the map $g$.  This completes the proof of the theorem for orthorhombic stiffness tensors.
\end{proof}

We note that the three solutions $(b_{12},b_{13}^*,b_{23}^*)$, $(b_{12}^*,b_{13}^*,b_{23})$, and $(b_{12}^*,b_{13},b_{23}^*)$ in~\eqref{eq:anomalouscompanions} are exactly the ``anomalous companions'' in~\cite[\S3]{HC2009}. Helbig and Carcione arrive at the existence of anomalous companions by making three quite reasonable assumptions that stiffness coefficients might satisfy in order for there to exist more than one set of stiffnesses that gives rise to the same slowness surface.  In other words, their conditions give sufficient conditions for the existence of anomalous companions. Our work shows that for a generic orthorhombic slowness polynomial, the anomalous companions in~\cite{HC2009} are the \emph{only} possible anomalous companions. \\

Before embraking on the proof of Theorem~\ref{thm:polynomial-to-tensor}(\ref{item:3Dmono}), we establish some notation around monoclinic stiffness tensors. Using Voigt notation~\eqref{eq:voigt3D}, such a tensor has 21 components $b_{ij}$, $1 \leq i \leq j \leq 6$, but
\[
b_{14} = b_{16} = b_{24} = b_{26} = b_{34} = b_{36} = b_{45} = b_{56} = 0,
\]
leaving at most~$13$ independent components $b_{11}$, $b_{12}$, $b_{13}$, $b_{15}$, $b_{22}$, $b_{23}$, $b_{25}$, $b_{33}$, $b_{35}$, $b_{44}$, $b_{46}$, $b_{55}$, $b_{66}$. The homogenized slowness polynomial of such a tensor has the form
\begin{equation}
\label{eq:genericmonoclinic}
\begin{split}
\tilde P(\vect{p}) &= \det(\Gamma(\vect{p}) - p_0^2I_3) \\
&=
    c_{1}p_1^6 + c_{2}p_1^5p_3 + c_{3}p_1^4p_2^2 + c_{4}p_1^4p_3^2 + c_{5}p_1^4p_0^2 + c_{6}p_1^3p_2^2p_3 + c_{7}p_1^3p_3^3 \\ &\quad+ c_{8}p_1^3p_3p_0^2
    + c_{9}p_1^2p_2^4 + c_{10}p_1^2p_2^2p_3^2 + c_{11}p_1^2p_2^2p_0^2 + c_{12}p_1^2p_3^4 + c_{13}p_1^2p_3^2p_0^2\\ &\quad+ c_{14}p_1^2p_0^4 + c_{15}p_1p_2^4p_3 + c_{16}p_1p_2^2p_3^3 + c_{17}p_1p_2^2p_3p_0^2 + c_{18}p_1p_3^5 + c_{19}p_1p_3^3p_0^2 \\ &\quad+ c_{20}p_1p_3p_0^4 + c_{21}p_2^6 + c_{22}p_2^4p_3^2 + c_{23}p_2^4p_0^2 + c_{24}p_2^2p_3^4 + c_{25}p_2^2p_3^2p_0^2\\ &\quad + c_{26}p_2^2p_0^4 + c_{27}p_3^6 + c_{28}p_3^4p_0^2 + c_{29}p_3^2p_0^4 + c_{30}p_0^6  
\end{split}    
\end{equation}
where, for example, we have
\begin{equation}
\label{eq:coeffexamplemonoclinic}
    c_{8} = -2b_{11}b_{35} - 2b_{11}b_{46} + 2b_{13}b_{15} - 2b_{15}b_{66} - 2b_{35}b_{66} - 2b_{46}b_{55}.
\end{equation}

\begin{proof}[Proof of Theorem~\ref{thm:polynomial-to-tensor}(\ref{item:3Dmono})] The proof is the same, \emph{mutatis mutandis} as the proof of Theorem~\ref{thm:polynomial-to-tensor}(\ref{item:3Dortho}). We point out the differences. 

The 13 stiffness parameters of a monoclinic tensor produce the 30 coefficients $c_1,\dots,c_{30}$ of the slowness polynomial~\eqref{eq:genericmonoclinic}, although $c_{30} = -1$, so we ignore it. Homogenizing with respect to the auxiliary variable $r$ as in~\eqref{eq:2Dreconstruction}, we arrive at a rational map
\[
\PP^{13}_{(b_{11}:\,\cdots\,:b_{66}:r)} \dasharrow \PP^{28}_{(\tilde c_1:\,\cdots\,:\tilde c_{29})}.
\]
We reduce the complexity of this rational map by leveraging the relations
\begin{equation}
\label{eq:linearrelationsmonoclinic}
\begin{split}
    c_{14} &= b_{11} + b_{55} + b_{66}, \\
    c_{20} &= b_{15} + b_{35} + b_{46}, \\
    c_{26} &= b_{22} + b_{44} + b_{66}, \\
    c_{29} &= b_{33} + b_{44} + b_{55}.
\end{split}    
\end{equation}
These relations imply that for a general set of stiffness parameters, the quantities $b_{22}$, $b_{44}$, $b_{46}$ and $b_{66}$ are determined uniquely by the values of $b_{11}$, $b_{15}$, $b_{33}$, $b_{35}$, $b_{55}$, $c_{14}$, $c_{20}$, $c_{26}$ and $c_{29}$. Accordingly, we simplify the map $\PP^{13} \dasharrow \PP^{28}$ by eliminating $b_{22}$, $b_{44}$, $b_{46}$, and $b_{66}$ to the map
\begin{align*}
g \colon\PP^9 &\dasharrow \PP^{24} \\
(b_{11}:b_{12}:b_{13}:b_{15}:b_{23}:b_{25}:b_{33}:b_{35}:b_{55}:r) &\mapsto \left(\tilde c_1:\,\cdots\,:\widehat{\tilde c}_{14}:\,\cdots\,:\widehat{\tilde c}_{20}:\,\cdots\,:\widehat{\tilde c}_{26}:\,\cdots\,:\tilde c_{28}:\widehat{\tilde c}_{29}\right).
\end{align*}
Now proceed as in the proof of Theorem~\ref{thm:polynomial-to-tensor}(\ref{item:3Dortho}): resolve the indeterminacy locus $\Pi$ of~$g$ through a blow-up process to get a surjective proper morphism $h\colon X \to Y$, and show there is a Zariski open subset of~$Y$ over which all fibers consist of \emph{at most two points}. For this last part, consider the point
\[
v = (b_{11}:b_{12}:b_{13}:b_{15}:b_{23}:b_{25}:b_{33}:b_{35}:b_{55}:r) =
(236:42:25:7:40:2:68:9:14:1) \in \PP^{9}(\Q),
\]
coming from the stiffness parameters for muscovite~\cite{muscovite}. We compute
\begin{align*}
y :&= h(u) = g(v) \\
&= ( 5853: 4746: -34628: -9848: 12102: 27175: -99912: 59245: 89980: -71422: -69036: \\
&\qquad 1726252: -1451432: -7151272: 864904: 6126476: 6721008: -813750: -9070032: \\
&\qquad-920304: 2329848: -1622864: -1078660: 8663468: -4103896 ) \in \PP^{24}(\Q).
\end{align*}
A Gr\"obner basis calculation included in~\cite{this-paper} shows that $h^{-1}(y)\subset \PP^9(\Q)$ consists of two closed points, namely
\begin{align*}
(236:42:25:7:40:2:68:9:14:1)&\quad\text{and}\quad (236:458:25:7:124:30:68:9:14:1),
\end{align*}
which includes the original stiffness parameters and one anomalous companion.

Next, we show that a non-empty fiber of the map $g$ has \emph{at least} two points. The above calculation suggests that anomalous companions can differ in the stiffnesses $b_{12}$, $b_{23}$, and $b_{25}$. Inspecting the coefficients of~\eqref{eq:genericmonoclinic} we see that
\begin{align*}
c_{11} &= b_{12}^2 + 2b_{12}b_{66} -b_{11}b_{22} - b_{11}b_{44}  + 2b_{15}b_{46} - b_{22}b_{55} + b_{25}^2 + 2b_{25}b_{46} - b_{44}b_{66} + b_{46}^2 - b_{55}b_{66}, \\
\label{eq:c8}
c_{25} &= b_{23}^2 + 2b_{23}b_{44} -b_{22}b_{33} - b_{22}b_{55}   + b_{25}^2 + 2b_{25}b_{46} - b_{33}b_{66} + 2b_{35}b_{46} - b_{44}b_{55} - b_{44}b_{66} + b_{46}^2.
\end{align*}
Considering these equations in turn as quadratics in $b_{12}$, $b_{23}$, and $b_{25}$ suggests that each of these stiffnesses can take on two possible values $b_{ij}$ and $b_{ij}^*$, related by the equations
\begin{align*}
b_{12} + b_{12}^* &= -2b_{66}, \\ 
b_{23} + b_{23}^* &= -2b_{44}, \\ 
b_{25} + b_{25}^* &= -2b_{46}.
\end{align*}
These equations are coupled by the linear relations~\eqref{eq:linearrelationsmonoclinic}, so only two possible triples can occur as stiffness parameters:
\begin{equation}
\label{eq:anomalouscompanionsmonoclinic}
(b_{12},b_{13},b_{23})\quad\text{and}\quad
(b_{12}^*,b_{13}^*,b_{23}^*)    
\end{equation}
Now we can simply verify that the two points
\begin{align*}
&(b_{11}:b_{12}:b_{13}:b_{15}:b_{23}:b_{25}:b_{33}:b_{35}:b_{55}:r),\quad\text{and}\\
&(b_{11}:-b_{12}-2b_{66}:b_{13}:b_{15}:-b_{23}-2b_{44}:-b_{25}-2b_{46}:b_{33}:b_{35}:b_{55}:r),
\end{align*}
have the same image under the map $g$.  This completes the proof of the theorem.
\end{proof}

%%%%%%%%%%%
\subsection{Which polynomials are slowness polynomials?}
\label{sec:slownesspolys}
In dimension~$2$, we have seen~\eqref{eq:gen2Dpoly} that the slowness polynomial has the form
\begin{equation}
\label{eq:general2Dpoly}
c_1p_1^4 + c_2p_1^3p_2 + c_3p_1^2p_2^2 + c_4p_1^2 + c_5p_1p_2^3 + c_6p_1p_2 + c_7p_2^4 + c_8p_2^2 + c_9
\end{equation}
for some $(c_1,\dots,c_9) \in \R^9$. However, not every polynomial of this kind arises from a stiffness tensor. For example, a close inspection of~\eqref{eq:gen2Dpoly} shows that we must have $c_9 = 1$.  Furthermore, the remaining coefficients $c_1,\dots,c_8$ are subject to the relations~\eqref{eq:2Dreconstructiondehom}. We can use elimination theory to compute the exact set of constraints that must be satisfied by $c_1,\dots,c_8$ (implicitly, from now on we simply take for granted that $c_9 = 1$). As a by-product, we shall obtain a third proof of Theorem~\ref{thm:polynomial-to-tensor}(\ref{item:2D}).  While in principle a similar argument could be used in the case $n = 3$, the required computations are currently infeasible.

Let~$X$ be the variety in the affine space $\Aff_\C^{14}$ with coordinates $b_{11},\dots,b_{33},c_1,\dots,c_8$ cut out by the equations~\eqref{eq:2Dreconstructiondehom}. More precisely, $X$ is $\Spec \C[b_{11},\dots,b_{33},c_1,\dots,c_8]/I$, where~$I$ is the ideal of $\C[b_{11},\dots,b_{33},c_1,\dots,c_8]$ given by
\[
\begin{split}
    I &\coloneqq \langle c_1 - (b_{11}b_{33}-b_{13}^2), c_2 - 2(b_{11}b_{23}-b_{12}b_{13}), \\
    &\qquad c_3 - (b_{11}b_{22}-b_{12}^2-2b_{12}b_{33}+2b_{13}b_{23}), c_4 + (b_{11}+b_{33}), \\
    &\qquad c_5 - 2(-b_{12}b_{23} + b_{13}b_{22}), c_6 + 2(b_{13}+b_{23}), \\
    &\qquad c_7 - (b_{22}b_{33} - b_{23}^2),  c_8 + (b_{22}+b_{33}) \rangle.
\end{split}
\]
We consider the two projections
\[
\xymatrix{
\Aff_\C^{14} = \Spec \C[b_{11},\dots,b_{33},c_1,\dots,c_8] \ar[rr]^q \ar[d]^p  & & \Aff_\C^6 = \Spec \C[b_{11},\dots,b_{33}] \\
\Aff_\C^{8} = \Spec \C[c_1,\dots,c_8]
}
\]
and by a slight abuse of notation, we also denote their restrictions to~$X$ by $p\colon X \to \Aff_\C^{8}$ and $q\colon X\to \Aff_\C^{6}$.  An elementary but important observation is that~$q\colon X \to \Aff_\C^6$ is an isomorphism, because the ring map
\[
\C[b_{11},\dots,b_{33},c_1,\dots,c_8] \to \C[b_{11},\dots,b_{33}]
\]
that sends~$b_{ij}$ to itself and maps~$c_i$ according to the relations~\eqref{eq:2Dreconstructiondehom} (so, e.g.,~$c_1$ maps to $b_{11}b_{33}-b_{13}^2$) is surjective and has kernel~$I$.  This tells us that~$X$ is a $6$-dimensional complex algebraic variety. 

We now turn to the projection $p\colon X \to  \Aff_\C^{8}$. The image $p(X)$ consists of $8$-tuples $(c_1,\dots,c_8)$ that, together with $c_9 = 1$, give a set of coefficients of a polynomial that is the slowness polynomial of at least one stiffness tensor (not necessarily positive) in dimension $2$. By~\cite[\S4.4, Theorem~4]{CLO}, the Zariski \emph{closure} of the image of~$p$ is cut out by the elimination ideal
\[
J \coloneqq I \cap \C[c_1,\dots,c_8] \subseteq \C[c_1,\dots,c_8],
\]
and a basis for this ideal can be extracted from an appropriate Gr\"obner basis for~$I$ by elimination theory (e.g.,~\cite[\S3.1, Theorem~2]{CLO}).  A {\tt magma} calculation~\cite{this-paper} shows that 
\begingroup
\allowdisplaybreaks
\begin{align*}
J = \langle 	
& -16c_1^2c_3 + 4c_1c_2^2 - 8c_1c_2c_5 + 16c_1c_3c_4c_8 - 4c_1c_3c_6^2 + 32c_1c_3c_7 - 16c_1c_3c_8^2 - 12c_1c_5^2 \\
			&\quad\quad	+ 16c_1c_5c_6c_8 - 16c_1c_6^2c_7 - 4c_2^2c_4c_8 + c_2^2c_6^2 - 12c_2^2c_7 + 4c_2^2c_8^2 + 16c_2c_4c_6c_7 \\
			&\quad\quad	- 2c_2c_5c_6^2 - 8c_2c_5c_7 - 4c_3c_4^2c_7 + 4c_3c_4c_7c_8 - c_3c_6^2c_7 - 4c_3c_7^2 + c_4^2c_5^2 - c_4c_5^2c_8 \\
			&\quad\quad	+ c_5^2c_6^2 + 4c_5^2c_7,\\
			& -4c_1^2 + 4c_1c_4c_8 - c_1c_6^2 + 8c_1c_7 - 4c_1c_8^2 - c_2^2 - 2c_2c_5 + 2c_2c_6c_8 - c_3c_6^2 - 4c_4^2c_7 \\
			&\quad\quad	+ 2c_4c_5c_6 + 4c_4c_7c_8 - c_5^2 - c_6^2c_7 - 4c_7^2 \rangle.
\end{align*}
\endgroup
With such an explicit description of~$J$, it is possible to compute the dimension of $Y \coloneqq \overline{p(X)} = \Spec \C[c_1,\dots,c_8]/J$. A {\tt magma} computation shows that the dimension is~$6$, which is the same dimension of~$X$. 

One can go further and compute the image $p(X)$, and not simply its Zariski closure, using an effective version of Chevalley's Theorem, which asserts that the set-theoretic image $p(X)$ is a constructible set~\cite{Barakat}. This way we obtain necessary and sufficient conditions on $(c_1,\dots,c_8)$ so that \eqref{eq:general2Dpoly} is the slowness polynomial for a stiffness tensor (note, however, that our algebro-geometric set-up does not take into account the positivity condition that must be satisfied by a physical stiffness tensor). For an ideal $I' \subset \C[c_1,\dots,c_8]$, write $V(I')$ for the affine variety cut out in $\Aff^8_\C$ by the ideal $I'$. Then, using the package {\tt ZariskiFrames}~\cite{ZariskiFrames}, we compute that
\[
p(X) = \left(V(I_1)\setminus V(J_1)\right)\ \cup\ 
\left(V(I_2)\setminus V(J_2)\right)\ \cup\ 
\left(V(I_3)\setminus V(J_3)\right),
\]
where
\begingroup
\allowdisplaybreaks
\begin{align*}
    I_1 &= J, \\
    J_1 &= \langle
        c_4^2 - 2c_4c_8 + c_6^2 + c_8^2,
    -c_2c_6 + 2c_3c_4 - 2c_3c_8 - 4c_4c_7 + 3c_5c_6 + 4c_7c_8,\\
    &\quad\quad c_2^2 - 6c_2c_5 + 4c_3^2 - 16c_3c_7 + 9c_5^2 + 16c_7^2,
    -3c_1c_6 + 2c_2c_4 - 2c_2c_8 + c_3c_6 + c_6c_7,\\
    &\quad\quad-c_1c_6 - c_3c_6 + 2c_4c_5 - 2c_5c_8 + 3c_6c_7,
    2c_1c_4 - 2c_1c_8 + c_2c_6 - 2c_4c_7 + c_5c_6 + 2c_7c_8,\\
    &\quad\quad c_1c_3 - 2c_1c_7 - c_2c_5 + c_3^2 - 5c_3c_7 + 3c_5^2 + 6c_7^2,\\
    &\quad\quad c_1c_2 - 3c_1c_5 + c_2c_3 - 3c_2c_7 + c_3c_5 + c_5c_7, \\
    &\quad\quad c_1^2 - 2c_1c_7 + 2c_2c_5 - c_3^2 + 4c_3c_7 - 2c_5^2 - 3c_7^2
    \rangle,\\
    I_2 &= \langle
        c_4^2 - 2c_4c_8 + c_6^2 + c_8^2,
    -c_1c_6 - c_3c_6 + 2c_4c_5 - 2c_5c_8 + 3c_6c_7,\\
    &\quad\quad -c_2c_6 + 2c_3c_4 - 2c_3c_8 - 4c_4c_7 + 3c_5c_6 + 4c_7c_8, \\
    &\quad\quad-3c_1c_6 + 2c_2c_4 - 2c_2c_8 + c_3c_6 + c_6c_7,
    2c_1c_4 - 2c_1c_8 + c_2c_6 - 2c_4c_7 + c_5c_6 + 2c_7c_8,\\
    &\quad\quad c_1c_6^2 - 16c_2c_5 + 4c_2c_6c_8 + 8c_3^2 - c_3c_6^2 - 32c_3c_7 + 16c_5^2 + 
        4c_5c_6c_8 - 7c_6^2c_7 + 32c_7^2,\\
    &\quad\quad 16c_1c_5 - 8c_1c_6c_8 - 4c_2c_3 + c_2c_6^2 + 8c_2c_7 - 4c_3c_5 + 8c_4c_6c_7 -
        c_5c_6^2 - 8c_5c_7,\\
    &\quad\quad c_1c_3 - 2c_1c_7 - c_2c_5 + c_3^2 - 5c_3c_7 + 3c_5^2 + 6c_7^2, \\
    &\quad\quad c_2^2 - 6c_2c_5 + 4c_3^2 - 16c_3c_7 + 9c_5^2 + 16c_7^2, \\
    &\quad\quad c_1c_2 - 3c_1c_5 + c_2c_3 - 3c_2c_7 + c_3c_5 + c_5c_7, \\
    &\quad\quad c_1^2 - 2c_1c_7 + 2c_2c_5 - c_3^2 + 4c_3c_7 - 2c_5^2 - 3c_7^2 \rangle, \\
J_2 &= \langle 
    c_4^2 - 2c_4c_8 + c_6^2 + c_8^2,
    4c_3 - 2c_4c_8 - c_6^2 + 24c_7 - 6c_8^2,\\
    &\quad\quad
    2c_2 - c_4c_6 + 2c_5 - c_6c_8,
    4c_1 - 2c_4c_8 + c_6^2 - 4c_7 + 2c_8^2,\\
    &\quad\quad 8c_4c_7 - 2c_4c_8^2 - 2c_5c_6 + c_6^2c_8 - 8c_7c_8 + 2c_8^3,\\
    &\quad\quad 2c_4c_5 - c_4c_6c_8 - 2c_5c_8 + 8c_6c_7 - c_6c_8^2,\\
    &\quad\quad 4c_5^2 - 4c_5c_6c_8 + c_6^2c_8^2 + 64c_7^2 - 32c_7c_8^2 + 4c_8^4
    \rangle,\\
    I_3 &= \langle
    c_8^2-4c_7,c_6c_8-2c_5,c_4c_8-c_1-c_3-c_7,2c_6c_7-c_5c_8,\\
    &\quad\quad c_6^2+2c_1-2c_3+2c_7,c_4c_6-2c_2,c_4^2-4c_1
    \rangle, \\
    J_3 &= \langle 1\rangle.
\end{align*}
\endgroup

Note that $V(J_3) = \emptyset$. 

\begin{proof}[Third proof of Theorem~\ref{thm:polynomial-to-tensor}(\ref{item:2D})] Since $p\colon X \to Y$ is a dominant morphism of integral schemes of finite type over a field, both of the same dimension, Chevalley's theorem~\cite[Exercise~II.3.22{e}]{Hartshorne} implies that there is a Zariski open subset $U \subset Y$ such that the fiber $p^{-1}(u)$ for $u \in U$ is a finite set.  In other words, for each $u \in U$, there are only finitely many possible values of $b_{11},\dots,b_{33}$ such that the relations~\eqref{eq:2Dreconstructiondehom} hold; more plainly, there are only finitely many stiffness tensors associated to a slowness polynomial corresponding to a point $u \in U$.  It is possible to choose~$U$ so that the number of stiffness tensors is \emph{constant} as one varies $u \in U$.  This constant is the degree of the map~$p$, which is equal to the degree of the function field extension $[\C(X):\C(Y)]$.  We use {\tt magma} to compute this quantity and show that it is~$1$; see~\cite{this-paper}.  The computation in fact gives explicit expressions for $b_{11},\dots,b_{33}$ in terms of $c_1,\dots,c_8$.  It shows that the map $p\colon X \to Y$ is a surjective, birational morphism, i.e.,~$p$ has an inverse defined on a Zariski open subset of~$Y$.
\end{proof}

Conjecture~\ref{conj:uniquereconstruction} can in principle be proved using the same template as above. However, the symbolic computations required when computing Gr\"obner bases are well beyond the capabilities of modern-day desktop computers. The slowness polynomials involved have $50$ monomials, with coefficients $c_1,\dots,c_{50}$, and the stiffness tensor has $21$ components $b_{11},\dots,b_{66}$. The analogous correspondence diagram for $n = 3$ has the form
\[
\xymatrix{
\Aff_\C^{71} = \Spec \C[b_{11},\dots,b_{66},c_1,\dots,c_{50}] \ar[d]^p \ar[rr]^q & & \Aff_\C^{21} = \Spec \C[b_{11},\dots,b_{66}]  \\
 \Aff_\C^{50} = \Spec \C[c_1,\dots,c_{50}] & &
}
\]
Using the map~$q$ as before we can show that the variety $X \subset \Aff_{\C}^{71}$ parametrizing slowness polynomials in terms of stiffness tensors has dimension~$21$. As before the closure $Y = \overline{p(X)}$ of the image of~$p$ could in principle be computed using elimination theory. This would give a set of polynomials generating an ideal~$J$ describing the closure of the image $\overline{p(X)}$.

%%%%%%%%%%%
\subsection{Positivity of anomalous companions for orthorhombic materials}
\label{sec:Cayley-cubic}
\begin{figure}
    \centering
    \includegraphics[clip, trim=3cm 3cm 3cm 3cm, width=.4\textwidth]{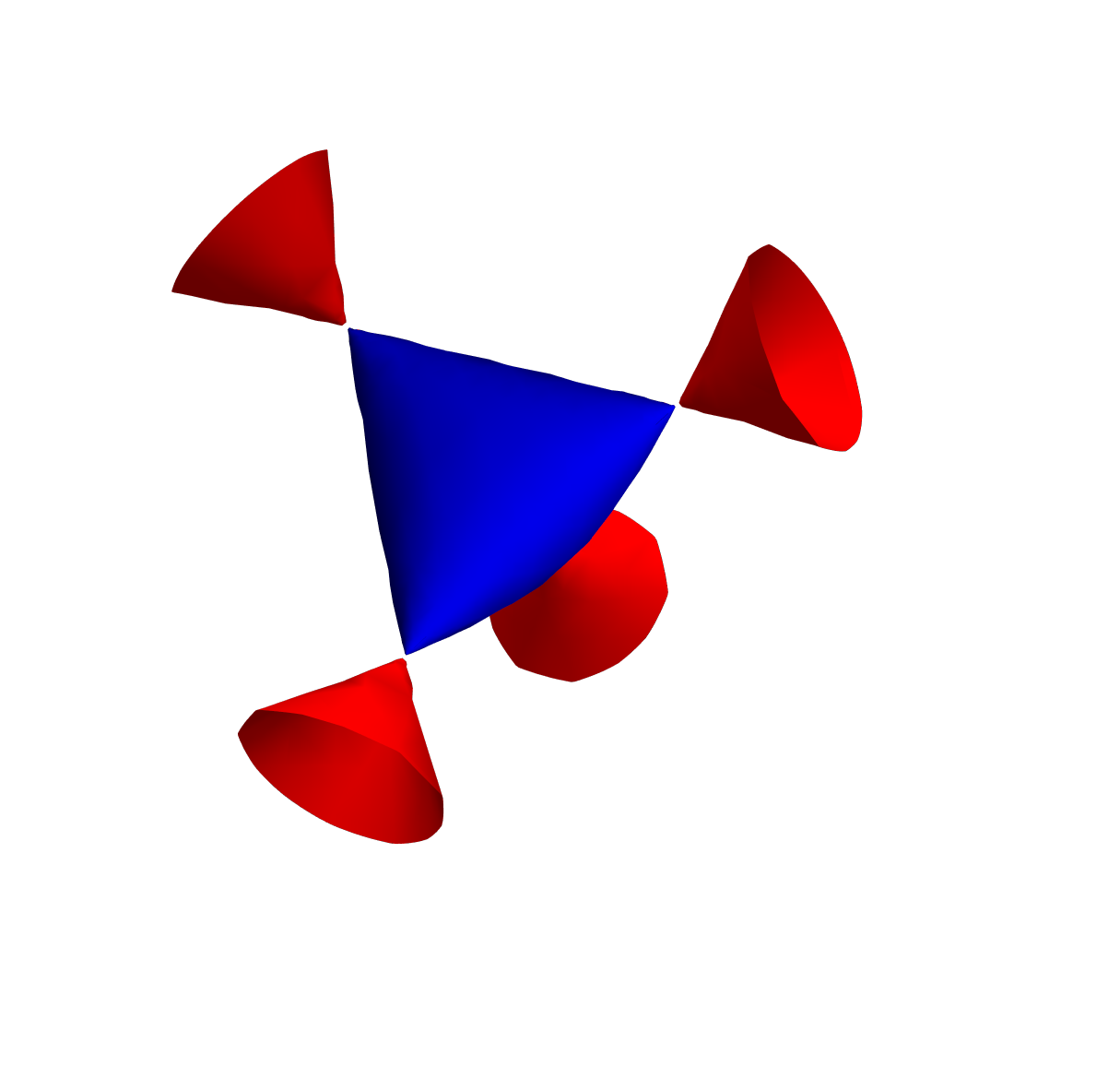}
    \caption{The Cayley cubic surface $1 + 2xyz - x^2 - y^2 - z^2 = 0$ in~$\R^3$.
    Points in the \textcolor{blue}{roughly tetrahedral shape in the middle} correspond to anomalous companions of an orthorhombic stiffness tensor that are physically realizable, while points in the \textcolor{red}{four cone-shaped regions extending to infinity} do not.}
    \label{fig:cayley}
\end{figure}

For completeness, we summarize here the analysis in~\cite{HC2009} characterizing which triples of stiffnesses~\eqref{eq:anomalouscompanions} give rise to positive orthorhombic stiffness tensors. Positivity requires that the $6\times 6$ matrix of stiffnesses
\[
\begin{pmatrix}
b_{11} & b_{12} & b_{13} & & & \\
b_{12} & b_{22} & b_{23} & & & \\
b_{13} & b_{23} & b_{33} & & & \\
 & & & b_{44} & & \\
 & & & & b_{55} & \\
 & & & & & b_{66} 
\end{pmatrix}
\]
be positive definite. By Sylvester's criterion~\cite[\S7.6]{Meyer}, this implies that
\[
b_{11} > 0, \quad b_{22} > 0,\quad b_{33} > 0,\quad b_{44} > 0,\quad b_{55} > 0,\quad b_{66} > 0,
\]
and that the $2\times 2$ minors
\[
b_{11}b_{22} - b_{12}^2,\quad b_{11}b_{33} - b_{13}^2,\quad b_{22}b_{33} - b_{23}^2
\]
are also positive, implying the inequalities
\begin{equation}
\label{eq:2minors}
-\sqrt{b_{11}b_{22}} < b_{12} < \sqrt{b_{11}b_{22}},\quad
-\sqrt{b_{11}b_{33}} < b_{13} < \sqrt{b_{11}b_{33}},\quad
-\sqrt{b_{22}b_{33}} < b_{23} < \sqrt{b_{22}b_{33}}.
\end{equation}
In addition, the $3\times 3$ leading principal minor must also be positive:
\begin{equation}
\label{eq:3minor}
    b_{11}b_{22}b_{33} + 2b_{12}b_{13}b_{23} - b_{11}b_{23}^2 - b_{22}b_{13}^2 - b_{33}b_{12}^2 > 0.
\end{equation}
Let 
\[
x \coloneqq \frac{b_{12}}{\sqrt{b_{11}b_{22}}},\quad
y \coloneqq \frac{b_{13}}{\sqrt{b_{11}b_{33}}},\quad\text{and}\quad
z \coloneqq \frac{b_{23}}{\sqrt{b_{22}b_{33}}}
\]
Then the conditions~\eqref{eq:2minors} and~\eqref{eq:3minor} become, respectively,
\begin{equation*}
    -1 < x < 1\quad -1 < y < 1, \quad -1 <z < 1
\end{equation*}
and
\begin{equation*}
    1 + 2xyz - x^2 - y^2 - z^2 > 0.
\end{equation*}
The affine surface $1 + 2xyz - x^2 - y^2 - z^2 = 0$ is the ubiquitous Cayley cubic surface! Positivity of an anomalous companion is equivalent to having the point corresponding to the companion lying inside the finite ``tetrahedral'' region in~$\R^3$ determined by the four singularities of the cubic surface. See Figure~\ref{fig:cayley}.

%%%%%%%%%%%%%%%%%%%%%%%%%%%%%%%%%%%%%%%%%%%%%%%%%%%%
%%%%%%%%%%%%%%%%%%%%%%%%%%%%%%%%%%%%%%%%%%%%%%%%%%%%
%%%%%%%%%%%%%%%%%%                                           %%%%%%%%%%%%%%%%%
%%%%%%%%%%%%%%%%%%		Bibliography		%%%%%%%%%%%%%%%%%
%%%%%%%%%%%%%%%%%%                                           %%%%%%%%%%%%%%%%%
%%%%%%%%%%%%%%%%%%%%%%%%%%%%%%%%%%%%%%%%%%%%%%%%%%%%
%%%%%%%%%%%%%%%%%%%%%%%%%%%%%%%%%%%%%%%%%%%%%%%%%%%%

\end{document}